\documentclass[11pt]{article}
\usepackage{graphicx}
\usepackage[a4paper,pdftex=true]{geometry}
\usepackage[latin1]{inputenc}
\usepackage[T1]{fontenc}
\usepackage{aeguill}
\usepackage{amssymb}
\usepackage{amssymb,array}
\usepackage{syntax}
\usepackage{mdwtab}
\usepackage{mathenv}
\usepackage{float}
\usepackage{lscape}
\usepackage{here}
\usepackage[a4paper]{geometry}
\usepackage{amssymb}
\usepackage{bm}
\usepackage{fancyhdr}
\usepackage{graphicx}
\usepackage{natbib}
\usepackage{fancyvrb}

\newtheorem{theorem}{Theorem}
\newtheorem{lemma}[theorem]{Lemma}
\newtheorem{corollary}[theorem]{Corollary}

\newtheorem{proposition}[theorem]{Proposition}

\newtheorem{remark}{Remark}

\newtheorem{definition}{Definition}

\newcommand{\Vect}[1]{
  \bm{#1}
}

\newcommand{\RR}{ \mathbb{R}}

\newcommand{\Given}[2]{
#1 | #2
}

\newcommand{\ent}[1]{
\ensuremath{H_{_{#1}}}
}

\newcommand{\IM}[2]{
\ensuremath{I_{_{#1, #2}}}
}

\newcommand{\dE}[2]{
\ensuremath{d_{_{#1 , #2}}^E}%
}

\newcommand{\dEempty}{
\ensuremath{d^E}%
}

\newcommand{\hE}[2]{
\NIB[D]{#1}{#2}
}

\newcommand{\HE}[2]{
\IB[D]{#1}{#2}
}

\newcommand{\hS}[2]{
\NIB[S]{#1}{#2}
}

\newcommand{\HS}[2]{
\DIB[S]{#1}{#2}
}

\newcommand{\NIB}[3][]{
\ensuremath{\delta^{#1}_{_{#2 , #3}}}
}
\newcommand{\dIB}[3][]{
\ensuremath{d^{#1}_{_{#2 , #3}}}
}
\newcommand{\NIBempty}[1][]{
\ensuremath{\delta^{#1}}%
}
\newcommand{\dIBempty}[1][]{
\ensuremath{d^{#1}}%
}

\newcommand{\NIBk}[3]{
\ensuremath{\delta^{(#1)}_{_{#2, #3}}}
}
\newcommand{\NIBkempty}[1]{
\ensuremath{\delta^{(#1)}}%
}

\newcommand{\IB}[3][]{
\ensuremath{\Delta^{#1}_{_{#2 , #3}}}
}
\newcommand{\DIB}[3][]{
\ensuremath{D^{#1}_{_{#2 , #3}}}
}
\newcommand{\IBempty}[1][]{
\ensuremath{\Delta^{#1}}%
}
\newcommand{\DIBempty}[1][]{
\ensuremath{D^{#1}}%
}

\newcommand{\IBk}[3]{
\ensuremath{\Delta^{(#1)}_{_{#2, #3}}}
}

\newcommand{\IBkempty}[1]{
\ensuremath{\Delta^{(#1)}}%
}

\newcommand{\pen}[3][]{
\ensuremath{C^{#1}_{_{#2 , #3}}}
}
\newcommand{\penk}[3]{
\ensuremath{C^{(#1)}_{_{#2 , #3}}}
}
\newcommand{\penkempty}[1]{
\ensuremath{C^{(#1)}}
}
\newcommand{\penempty}[1][]{
\ensuremath{C^{#1}}
}

\newcommand{\dI}[2]{
\ensuremath{d_{_{#1 , #2}}^{^I}}
}
\newcommand{\dIempty}{
\ensuremath{d^{^I}}%
}

\newcommand{\DE}[2]{
\ensuremath{D_{_{#1 , #2}}^{^E}}
}

\newcommand{\DEempty}{
\ensuremath{D^{^E}}%
}

\newcommand{\DI}[2]{
\ensuremath{D_{_{#1 , #2}}^{^I}}
}
\newcommand{\DIempty}{
\ensuremath{D^{^I}}%
}

\newcommand{\nibd}[1][ ]{NIB-divergence#1}

\newcommand{\ibDiv}[1][ ]{information-based divergence#1}
\newcommand{\ibd}[1][ ]{IB-divergence#1}

\newcommand{\compl}[1][ ]{complexity term#1}
\newcommand{\Compl}[1][ ]{Complexity term#1}

\title{\Large \bf \sc  Normalized information-based divergences}

\author{ {\sc By J.-F. Coeurjolly, R. Drouilhet and J.-F. Robineau} \\
{\it University of Grenoble 2, France}}
\begin{document}
\maketitle


\vspace*{.5cm}

\noindent {\it Abstract}

\begin{center}
{\small \begin{minipage}{12cm}
This paper is devoted to the mathematical study of some divergences based on the mutual information well-suited to categorical random vectors. These divergences are generalizations of the ``entropy distance" and ``information distance". Their main characteristic is that they combine a complexity term and the mutual information. We then introduce the notion of (normalized) \ibDiv[,] propose several examples and discuss their mathematical properties in particular in some prediction framework.
\end{minipage}}
\end{center}

\vspace*{.5cm}

\noindent {\it Keywords and phrases}. Information theory, entropy distance, information distance, triangular inequality, redundancy.

\newpage

\section{Introduction}

Shannon information theory, usually just called {\it information} theory was introduced in 1948, \cite{Shannon48b}.
The theory aims at providing a means for measuring information. More precisely, the amount of information in an object may be measured by its \emph{entropy} and may be interpreted as the length of the description of the object by some encoding way. In the Shannon approach, the objects to be encoded are assumed to be outcomes of a known source.
Shannon theory also provides the notion of {\it mutual information} (related to two objects) which plays a central role in many applications, from lossy compression to machine learning methods.

Several authors noticed that it would be useful to modify the mutual information such that the resulting quantity becomes a metric in a strict sense.
As a first example, \cite{crutchfield:90}, \cite{hillman-formal} introduced the {\it entropy distance} defined as the sum of the conditional entropies.
Other interesting measures are the {\it information distance} \cite{journals/tit/BennettGLVZ98} and its normalized version named {\it similarity metric} introduced by \cite{LiEtAl04c} in the context of the Kolmogorov complexity theory. More precisely, the information distance is defined as the maximum of the conditional Kolmogorov complexities. The similarity metric is universal in the sense defined by the authors and is not computable, since it is based on the uncomputable notion of Kolmogorov complexity.

Recent papers have demonstrated useful application of suitable version of the similarity metric in areas as diverse as genomics, virology, languages, literature, music, handwritten digits and astronomy, \cite{CilVit05}. To apply the metric to real data, the authors have to replace the use of the noncomputable Kolmogorov complexity by an approximation using standard real-world compressors~: GenCompress for genomics, \cite{LiEtAl01}, the \emph{Normalized Compression Distance (NCD)} for music clustering, \cite{journals/corr/cs-SD-0303025}, the \emph{Normalized Google Distance (NGD)} for automatic meaning discovery, \cite{oai:eprints.pascal-network.org:1826}, are examples of effective compressors. To include the information distance and the similarity metric in a framework based on information theory concepts, we make use of the principle that \emph{expected Kolmogorov complexity equals Shannon entropy} and interested reader can refer to \cite{journals/corr/cs-IT-0410002}, \cite{Leung:78}, \cite{HammerEtAl00} for more details. 
Consequently, the entropy and information distances are both expressed in terms of conditional entropies: the first one as their sum and the second one as their maximum.
\cite{Kraskov03} gives a proof of the triangular inequality for these distances and their respective normalized versions.

In the supervised learning framework, the use of some selection method of covariables amoung a large number is required when it is assumed that the data size is too small with respect to the number of the  available covariables in order to apply any existing discriminant analysis method. Such a problem has been widely treated, \cite{LiuMotoda98}. 
The approach undertaken by \cite{Robineau04} is mainly based on three kinds of methodological tools. The first one is a supervised quantization method consisting in the simplification of covariables too complex (in particular with a too large number of possible values). Indeed, our main belief is that, in order to predict the class variable generally representing a small number of categories of data, each possibly predictive covariable must not be too complex. The second one is a more usual step by step selection method  combining the simplified covariables together in order to detect cluster of data of the same class. The last one is aimed at detecting redundancy among the covariables set. These three tasks may be realized using the entropy or information distances (or their normalized versions). Let us emphasize some properties allowing to understand the usefulness of these criterions in such a context. 
The entropy and information distances $\DEempty$ and $\DIempty$ can be rewritten as the difference between some term, respectively the joint entropy and the maximum of the marginal entropies, and the mutual information. The first term may be interpreted as a complexity term.
Moreover, both are independence measures with the particular property to be minimal (in fact equal to 0) when random vectors share exactly the same information.
\cite{Robineau04} proposes then to extend the definition of the entropy and information distances by introducing  the notion of information-based divergence $\IB{\Vect{X}}{\Vect{Y}}$ between two categorical random vectors $\Vect{X}$ and  $\Vect{Y}$ defined as the difference of some complexity term $\pen{\Vect{X}}{\Vect{Y}}$ and the mutual information $\IM{\Vect{X}}{\Vect{Y}}$ and such that  $\pen{\Vect{X}}{\Vect{Y}}$ is an upper bound of $\IM{\Vect{X}}{\Vect{Y}}$ reached when $\Vect{X}$ and  $\Vect{Y}$ share exactly the same information. The notion of normalized information-based divergence $\NIB{\Vect{X}}{\Vect{Y}}$ derives directly by dividing the associated information-based divergence $\IB{\Vect{X}}{\Vect{Y}}$ by the complexity term $\pen{\Vect{X}}{\Vect{Y}}$. The normalized version $\dEempty$ and  $\dIempty$ of $\DEempty$ and $\DIempty$ are particular examples.
Other examples are given in \cite{Robineau04}.
Amoung them, one is of particular interest since its complexity term
$\penempty[S]$ is the mean of the marginal entropies. The associated (non-normalized) information-based divergence $\IBempty[S]$ is not so different from $\DEempty$ since it corresponds to its half. Nevertheless, the expression of its \compl $\penempty[S]$ really differs from the complexity term $\penempty[E]$  of $\DEempty$  (i.e. the joint entropy). For pratical purposes, we may argue that $\DIempty$, $\DEempty$ and $\IBempty[S]$ are not well-suited in prediction framework since a small value of these distances means that both the explained and explicative variables have a good knowledge of each other. This is  due to the fact that both conditional entropies have at least the same weight.

In this paper, this drawback is weakened by introducing a natural extension $\penempty[S,\alpha]$ of the \compl $\penempty[S]$ defined as a weighted  mean (by $\alpha$ and $(1-\alpha)$ for some $0<\alpha \leq 1$) of the minimum and maximum of marginal entropies. This kind of \compl leads to an expected \ibd[], \IBempty[S,\alpha] which is the weighted mean of the minimum and maximum of conditional entropies.

The paper is organized as follows. In Section~\ref{sec-dEdI}, we recall the definition and their main properties of the entropy and information distances (and their normalized version). Similarly to \cite{Granger04}, we extract the main characteristics to define some general concept of information divergence which could be theoretically applied in a more general setting (continuous, discrete, \ldots). Section~\ref{sec-ibd} concentrates itself on categorical data (and in particular discrete) random vectors, as it is usually the case in most of applications using entropy or information distance. We give the definition of (normalized) \ibDiv and propose several examples. We study their mathematical properties in a general context and propose some sufficient conditions for these divergences to verify some triangular's type inequality. Finally, in Section~\ref{sec-pred}, we exhibit some properties of information-based divergences in the special prediction framework. In particular, we show that these divergences are useful to detect redundancy.

\section{Normalized entropy distance and normalized information distance} \label{sec-dEdI}

Let us denote by $\Gamma$ the set of categorical random vectors, that is, discrete-valued random vectors with finite entropy. In the sequel, $\Vect{X},\Vect{Y}$ and $\Vect{Z}$ are three elements of such a set~$\Gamma$.

\subsection{Some notation}

We denote by $\ent{\Vect{X}}$ (when it exists) the Shannon entropy of $\Vect{X}$ given by 
$$
\ent{\Vect{X}}= - \sum_{\Vect{x} \in \Omega_{\Vect{X}}} p_{\Vect{X}}(\Vect{x}) \log( p_{\Vect{X}}(\Vect{x})) \qquad \mbox{ with } p_{\Vect{X}}(\Vect{x})=\mathbb{P}(\Vect{X}=\Vect{x}), 
$$
In the same way, one can define the joint entropy of $\Vect{X}$ and $\Vect{Y}$ denoted by $\ent{\Vect{X},\Vect{Y}}$, the conditional entropy of $\Vect{X}$ (resp. $\Vect{Y}$) by $\Vect{Y}$ (resp. $\Vect{X}$) denoted by $\ent{\Given{\Vect{X}}{\Vect{Y}}}$ (resp. $\ent{\Given{\Vect{Y}}{\Vect{X}}}$). Finally, we denote by $\IM{\Vect{X}}{\Vect{Y}}$ the mutual information between the random vectors $\Vect{X}$ and $\Vect{Y}$. When these different quantities exist, the following relations hold (see {\it e.g.} \cite{Cover91}):
\begin{eqnarray}
\ent{\Vect{X},\Vect{Y}} &=& \ent{\Vect{X}} + \ent{\Given{\Vect{Y}}{\Vect{X}}} = \ent{\Vect{Y}}+ \ent{\Given{\Vect{X}}{\Vect{Y}}} \label{eq-entJointe}\\
\IM{\Vect{X}}{\Vect{Y}} &=& \ent{\Vect{X}}-\ent{\Given{\Vect{X}}{\Vect{Y}}}= \ent{\Vect{Y}}-\ent{\Given{\Vect{Y}}{\Vect{X}}} = \ent{\Vect{X}}+\ent{\Vect{Y}}-\ent{\Vect{X},\Vect{Y}}\label{eq-defIM}
\end{eqnarray}

\subsection{Definition and characteristics}

We now shall present some measures allowing to overcome some drawbacks of the mutual information. As a first generalization, several authors noticed that it would be useful to modify the mutual information such that the resulting quantity becomes a metric in a strict sense. Two such measures exist and are well-known in the litterature. The first one called ``entropy distance'' is derived from the domain of information theory. The second one called ``information distance'' originates in works around the Kolmogorov complexity. Both measures are defined (when they exist) for two random vectors $\Vect{X}$ and $\Vect{Y}$ by:
\begin{list}{$\bullet$}{}
\item Entropy distance: 
\begin{equation} \label{DistE}
\DE{\Vect{X}}{\Vect{Y}} = \ent{\Given{\Vect{X}}{\Vect{Y}}}  + \ent{\Given{\Vect{Y}}{\Vect{X}}}
\end{equation}
\item Information distance: 
\begin{equation} \label{DistI}
\DI{\Vect{X}}{\Vect{Y}} = \max \left( \ent{\Given{\Vect{X}}{\Vect{Y}}}  , \ent{\Given{\Vect{Y}}{\Vect{X}}} \right).
\end{equation}
\end{list}
Both measures are indeed some modifications of mutual information since from~(\ref{eq-entJointe}) and~(\ref{eq-defIM}), we have
\begin{equation}\label{eq-DEDIIM}
\DE{\Vect{X}}{\Vect{Y}} = \ent{\Vect{X},\Vect{Y}} - \IM{\Vect{X}}{\Vect{Y}} \qquad \mbox{ and } \qquad 
\DI{\Vect{X}}{\Vect{Y}} = \max \left( \ent{\Vect{X}}, \ent{\Vect{Y}} \right) - \IM{\Vect{X}}{\Vect{Y}}.
\end{equation}
The quantities $\ent{\Vect{X},\Vect{Y}}$ and $\max\left( \ent{\Vect{X}},\ent{\Vect{Y}}\right)$ are upper-bounds of the mutual information $\IM{\Vect{X}}{\Vect{Y}}$ that are reached when $\Vect{X}$ and $\Vect{Y}$ share exactly the same information. In other words, these two measures are nonnegative and vanish if and only if $\ent{\Given{\Vect{Y}}{\Vect{X}}}= \ent{\Given{\Vect{X}}{\Vect{Y}}}=0$ expressing the fact that $\Vect{X}$ (resp. $\Vect{Y}$) predicts $\Vect{Y}$ (resp. $\Vect{X}$) with probability $1$. 

These measures satisfy 
\begin{equation} \label{eq-majoDist}
\DE{\Vect{X}}{\Vect{Y}} \leq \ent{\Vect{X},\Vect{Y}} \qquad \mbox{ and } \qquad \DI{\Vect{X}}{\Vect{Y}}\leq \max \left( \ent{\Vect{X}},\ent{\Vect{Y}} \right),
\end{equation}
where the equality holds if the vectors $\Vect{X}$ and $\Vect{Y}$ are independent. As noticed by~\cite{Kaltchenko04}, \cite{LiVit97} argued that in Bioinformatics an unnormalized distance may not be a proper evolutionary distance measure. It would put two long and complex sequences that differ only by a tiny fraction of the total information as dissimilar as two short sequences that differ by the same absolute amount and are completely random with respect to one another. 
To overcome this problem within the algorithmic framework \cite{LiVit97} form two normalized versions of distances $\DEempty$ and $\DIempty$. Their Shannon version have been proposed and studied by~\cite{Kraskov03}

\begin{definition} When they exist, one defines the two following measures:
\begin{list}{$\bullet$}{}
\item Normalized entropy distance:
$$ \dE{\Vect{X}}{\Vect{Y}} = \frac{  \ent{\Given{\Vect{X}}{\Vect{Y}}}+    \ent{\Given{\Vect{Y}}{\Vect{X}}}}{ \ent{\Vect{X},\Vect{Y}}}$$
\item Normalized information distance:
$$
\dI{\Vect{X}}{\Vect{Y}}= \frac{\max \left( \ent{\Given{\Vect{X}}{\Vect{Y}}}  , \ent{\Given{\Vect{Y}}{\Vect{X}}} \right)}{ \max \left(\ent{\Vect{X}}  , \ent{\Vect{Y}}\right) }
$$
\end{list}
Since $\ent{\Vect{X},\Vect{Y}}=0 \Leftrightarrow  \ent{\Vect{X}}= \ent{\Vect{Y}}=0 \Leftrightarrow \max \left(\ent{\Vect{X}}  , \ent{\Vect{Y}} \right)=0 $, we set by convention $ \dE{\Vect{X}}{\Vect{Y}}=0$ (resp. $ \dI{\Vect{X}}{\Vect{Y}}=0$) when $\ent{\Vect{X}}= \ent{\Vect{Y}}=0$.
\end{definition}

We are encouraged to define the following class of equivalence: the vectors $\Vect{X}$ and $\Vect{Y}$ are said to be equivalent if $\Vect{X}$ (resp. $\Vect{Y}$) predicts $\Vect{Y}$ (resp. $\Vect{X}$) with probability~1 and one will denote
\begin{equation} \label{eq-classeEqu}
\Vect{X} \sim \Vect{Y} \Leftrightarrow \ent{\Given{\Vect{Y}}{\Vect{X}}}= \ent{\Given{\Vect{X}}{\Vect{Y}}}=0 \Leftrightarrow \IM{\Vect{X}}{\Vect{Y}}=\ent{{\Vect{X},\Vect{Y}}}=\ent{\Vect{X}}=\ent{\Vect{Y}}
\end{equation}
Due to the previous convention
$$
\dE{\Vect{X}}{\Vect{Y}} = 0 \Leftrightarrow \dI{\Vect{X}}{\Vect{Y}}=0  \Leftrightarrow \Vect{X}\sim\Vect{Y}.
$$
From (\ref{eq-entJointe}) and~(\ref{eq-defIM}), one can obtain the following expressions for these two measures allowing some new interpretations.

\begin{proposition} \label{propdEdI}
We have the following expressions for $\dE{\Vect{X}}{\Vect{Y}}$ and $\dI{\Vect{X}}{\Vect{Y}}$.
\begin{eqnarray}
\dE{\Vect{X}}{\Vect{Y}} &=& 1 - \frac{\IM{\Vect{X}}{\Vect{Y}}}{\ent{{\Vect{X},\Vect{Y}}}} \label{expr1dE} \\\dI{\Vect{X}}{\Vect{Y}} &=& 1 - \frac{\IM{\Vect{X}}{\Vect{Y}}}{ \max \left(\ent{\Vect{X}},\ent{\Vect{Y}}\right)} \label{expr1dI}  \\&=& \max \left( 
\frac{\ent{\Given{\Vect{X}}{\Vect{Y}}} }{\ent{\Vect{X}}}, 
\frac{\ent{\Given{\Vect{Y}}{\Vect{X}}} }{\ent{\Vect{Y}}} \right) \label{expr3dI}
\end{eqnarray} 
\end{proposition}

\begin{proposition} \label{prop-dEdIDist}
The measures $\dEempty$ et $\dIempty$ constitute two distances bounded by 1.
\end{proposition}

To our knowledge, these results have been proved by \cite{Kraskov03}. Proofs are very similar to proofs of \cite{LiVitanyi03} who consider the algorithmic version of these distances. The proof is then omitted, but in Section~\ref{sec-infoMesDist}, we propose a result extending this one in the sense that we give conditions on measures that can be written as~(\ref{expr1dE}) and~(\ref{expr1dI}) to constitute a metric.

\subsection{Concept of information divergence}

We can exhibit from the previous study related to $\DIempty$, $\DEempty$, $\dIempty$ and $\dEempty$, some characteristics useful for an attempt to define the concept of information divergence denoted by $\IBempty$ in a more general setting. Let us first consider a similarity measure $\mathcal{I}_{\Vect{X},\Vect{Y}}$ (not necessarily the mutual information) minimal (in fact equal to 0) when $\Vect{X}$ and $\Vect{Y}$ are independent, and maximal (in fact equal to $\mathcal{I}_{\Vect{X},\Vect{X}}=\mathcal{I}_{\Vect{Y},\Vect{Y}}$) when the distributions of $\Vect{X}$ given $\Vect{Y}=\Vect{y}$ and $\Vect{Y}$ given $\Vect{X}=\Vect{x}$ are trivial. An information divergence $\IB{\Vect{X}}{\Vect{Y}}$ could satisfy the following properties:
\begin{list}{}{}
\item \textbf{[P1]} symmetry: $\IB{\Vect{X}}{\Vect{Y}}=\IB{\Vect{Y}}{\Vect{X}}$.
\item \textbf{[P2]} nonnegativeness: $\IB{\Vect{X}}{\Vect{Y}}\geq 0$.
\item \textbf{[P3]}  $\IB{\Vect{X}}{\Vect{Y}}$ is minimum (i.e. $\IB{\Vect{X}}{\Vect{Y}}= 0$) if and only if  $\Vect{X}$ and $\Vect{Y}$ share exactly the same information (i.e. $\mathcal{I}_{\Vect{X},\Vect{Y}}$ is maximal).
\item \textbf{[P4]}  $\IB{\Vect{X}}{\Vect{Y}}$ is maximum if and only if $\Vect{X}$ and $\Vect{Y}$ are independent (i.e.  $\mathcal{I}_{\Vect{X},\Vect{Y}}=0$).
\end{list}
Other supplementary properties could be that $\IB{\Vect{X}}{\Vect{Y}}$:
\begin{list}{}{}
\item \textbf{[P5]} is normalized: $\IB{\Vect{X}}{\Vect{Y}}\in[0,1]$ and $\IB{\Vect{X}}{\Vect{Y}}=1$ when $\Vect{X}$ and $\Vect{Y}$ are independent.
\item \textbf{[P6]} satisfies a triangular inequality: $\IB{\Vect{X}}{\Vect{Y}} \leq \IB{\Vect{X}}{\Vect{Z}} + \IB{\Vect{Z}}{\Vect{Y}}$.
\item\textbf{[P7]} invariant under continuous and strictly increasing transformations $\Vect{\varphi}(\cdot)$, $\Vect{\psi}(\cdot)$ of the vectors $\Vect{X}$ and $\Vect{Y}$, whenever they are quantitative random vectors.
\end{list}

There exists a large litterature on the discussion of criteria satisying the previous stated properties. We may cite \cite{Ullah96}, or a recent work of \cite{Granger04} who propose to detect the dependence between two possibly nonlinear processes through the Bhattacharya-Matusita-Hellinger measure of dependence given by
$$
S_{\rho} = \frac12 \int \int \left( \sqrt{f_1(\Vect{x},\Vect{y})} - \sqrt{f_2(\Vect{x},\Vect{y})} \right)^2 d\Vect{x}d\Vect{y},
$$
where $f_1$ (resp. $f_2$) is the joint density (resp. the product of marginal densities) of $\Vect{X}$ and $\Vect{Y}$. This measure, that has the other advantage to be applicable to continuous or discrete variables, satisfies properties \textbf{[P1]}-\textbf{[P7]} 
(in fact let us precise that \textbf{[P7]} 
is only valid if $\Vect{\varphi}(\cdot)=\Vect{\psi}(\cdot)$).

In some framework where the purpose is to predict some reference variable, one may find interesting to work with a divergence $\IB{\Vect{X}}{\Vect{Y}}$ which combines the minimization of a nonnegative complexity term denoted by $\mathcal{C}_{\Vect{X},\Vect{Y}}$ and the maximization of a nonnegative information term $\mathcal{I}_{\Vect{X},\Vect{Y}}$. The quantity $\mathcal{C}_{\Vect{X},\Vect{Y}}$ is called a complexity term since it is assumed to be expressed as a function of $\mathcal{H}_{\Vect{X}}$, $\mathcal{H}_{\Vect{Y}}$ and $\mathcal{H}_{\Vect{X},\Vect{Y}}$ measuring in some way respectively the complexity of vectors $\Vect{X}$, $\Vect{Y}$ and $({\Vect{X},\Vect{Y}})$. In other words, we may expect that an information divergence $\IB{\Vect{X}}{\Vect{Y}}$ could also satisfy the following properties:
\begin{list}{}{}
\item \textbf{[P8]} When $\Vect{X}_1$ and $\Vect{X}_2$ have the same complexity (in the sense that $\mathcal{C}_{\Vect{Y},\Vect{X}_1}= \mathcal{C}_{\Vect{Y},\Vect{X}_2}$): $\IB{\Vect{Y}}{\Vect{X}_1} < \IB{\Vect{Y}}{\Vect{X}_2}$ whenever $\Vect{X}_1$ has a better knowledge about $\Vect{Y}$ than $\Vect{X}_2$ (i.e. $\mathcal{I}_{\Vect{Y},\Vect{X}_1}> \mathcal{I}_{\Vect{Y},\Vect{X}_2}$).
\item \textbf{[P9]} When $\Vect{X}_1$ and $\Vect{X}_2$ have the same knowledge about $\Vect{Y}$ (i.e. $\mathcal{I}_{\Vect{Y},\Vect{X}_1}= \mathcal{I}_{\Vect{Y},\Vect{X}_2}$): $\IB{\Vect{Y}}{\Vect{X}_1}< \IB{\Vect{Y}}{\Vect{X}_2}$ whenever $\Vect{X}_1$ is simpler than $\Vect{X}_2$ in the sense that  $\mathcal{C}_{\Vect{Y},\Vect{X}_1}< \mathcal{C}_{\Vect{Y},\Vect{X}_2}$. Moreover, in this particular situation the fact that
\item[] \textbf{[P10]} $\mathcal{C}_{\Vect{Y},\Vect{X}_1}\leq \mathcal{C}_{\Vect{Y},\Vect{X}_2}$ must be equivalent to $\mathcal{H}_{\Vect{X}_1}\leq  \mathcal{H}_{\Vect{X}_2}$.
\item \textbf{[P11]} When $\Vect{X}_1$ and $\Vect{X}_2$ share almost exactly the same information (i.e. $\mathcal{I}_{\Vect{X}_1,\Vect{X}_2}$ is almost maximal and $\IB{\Vect{X}_1}{\Vect{X}_2}\simeq 0$) then the difference between the divergences $\IB{\Vect{Y}}{\Vect{X}_1}$ and $\IB{\Vect{Y}}{\Vect{X}_2}$ is almost zero (i.e. $\IB{\Vect{Y}}{\Vect{X}_1}\simeq \IB{\Vect{Y}}{\Vect{X}_2}$).
\end{list}
A class of candidates that satisfy \textbf{[P8]} and \textbf{[P9]} the previous statements could be of the form:
\begin{equation} \label{eq-defGenIB}
\IB{\Vect{X}}{\Vect{Y}} = \frac{\mathcal{C}_{\Vect{X},\Vect{Y}} -\mathcal{I}_{\Vect{X},\Vect{Y}}}{ \mathcal{W}_{\Vect{X},\Vect{Y}} },
\end{equation}
where $\mathcal{W}_{\Vect{X},\Vect{Y}}$ is a positive term. When $\mathcal{W}_{\Vect{X},\Vect{Y}}=\mathcal{C}_{\Vect{X},\Vect{Y}}$ we obtain a normalized information divergence. The properties \textbf{[P2]}-\textbf{[P3]} 
and the form~(\ref{eq-defGenIB}) implies that $\mathcal{C}_{\Vect{X},\Vect{Y}}$ is an upper bound of $\mathcal{I}_{\Vect{X},\Vect{Y}}$ reached when $\Vect{X}$ and $\Vect{Y}$ share exactly the same information.

In the rest of this paper we concentrate ourself on criteria described by~(\ref{eq-defGenIB}) that are in addition well-suited to categorical random variables (and in particular discrete random variables). In such a framework, we shall only describe some entropic-based criteria (i.e. $\mathcal{H}_{\Vect{X}}=\ent{\Vect{X}}$), and so the information term will be set to the mutual information $\IM{\Vect{X}}{\Vect{Y}}$.

\section{Information-based divergences and their normalized versions} \label{sec-ibd}

\subsection{Definition and examples}


\begin{definition} \label{def-infoMes}
Two criteria $\IBempty$ and $\NIBempty$ are respectively called an information-based divergence and a normalized information-based divergence (in short \ibd and \nibd) if they can respectively be written 
\begin{eqnarray}
 \IB{\Vect{X}}{\Vect{Y}} &=& \pen{\Vect{X}}{\Vect{Y}} -\IM{\Vect{X}}{\Vect{Y}} \label{eq-defIB} \\
\NIB{\Vect{X}}{\Vect{Y}} &=& \frac{ \pen{\Vect{X}}{\Vect{Y}} - \IM{\Vect{X}}{\Vect{Y}}}{\pen{\Vect{X}}{\Vect{Y}}} =1- \frac{\IM{\Vect{X}}{\Vect{Y}}}{\pen{\Vect{X}}{\Vect{Y}}}  \label{eq-defInfoMes} 
\end{eqnarray}
where the term $\pen{\Vect{X}}{\Vect{Y}}$ constitutes a \compl satisfying 
\begin{itemize}
\item[$(i)$] $\pen{\Vect{X}}{\Vect{Y}}=\pen{\Vect{Y}}{\Vect{X}}$
\item[$(ii)$] $\IM{\Vect{X}}{\Vect{Y}}\leq \pen{\Vect{X}}{\Vect{Y}}$ and this bound is achieved if and only if the random vectors $\Vect{X}$ and $\Vect{Y}$ are equivalent, i.e. if and only if $\Vect{X} \sim \Vect{Y}$.
\end{itemize}
We set by convention $\NIB{\Vect{X}}{\Vect{Y}}=0$ when $\pen{\Vect{X}}{\Vect{Y}}=\IM{\Vect{X}}{\Vect{Y}}=0$.

\end{definition}

This definition implies automatically that an \ibd $\IB{\Vect{X}}{\Vect{Y}}$ (resp. a \nibd $\NIB{\Vect{X}}{\Vect{Y}}$) satisfies properties \textbf{[P1]}-\textbf{[P4]} (resp. \textbf{[P1]}-\textbf{[P5]}). In the rest of the paper, the term $\pen{\Vect{X}}{\Vect{Y}}$ is expressed as

\begin{equation} \label{CxyfC}
\pen{\Vect{X}}{\Vect{Y}} =f_C\left( \ent{\Given{\Vect{X}}{\Vect{Y}}} ,\ent{\Given{\Vect{Y}}{\Vect{X}}},\IM{\Vect{X}}{\Vect{Y}} \right),
\end{equation}
where $f_C(\cdot,\cdot,\cdot)$ is a nonnegative function. Under such an expression of $\pen{\Vect{X}}{\Vect{Y}}$, the property \textbf{[P7]} is ensured since the conditional entropies and the mutual information depend only on the joint probability distribution of the categorical random vectors $\Vect{X}$ and $\Vect{Y}$.

From now on, we propose a series of examples for which we adopt the following convention: an \ibd (resp. a  \nibd) satisfying the triangular inequality is denoted \DIBempty (resp. \dIBempty) rather than \IBempty (resp. \NIBempty) . Moreover, each example will be particularized by some discriminating additonal letter in the same manner as $\DEempty$ and $\DIempty$ (resp. $\dEempty$ and $\dIempty$) which clearly constitute \ibd[s] (resp. \nibd[s]). 

In~\cite{Robineau04}, we investigate about two new entropic criteria naturally expressed by
$$
\hE{\Vect{X}}{\Vect{Y}} = \frac{1}{2} \left(\frac{\ent{\Given{\Vect{X}}{\Vect{Y}}}}{\ent{\Vect{X}}} + \frac{\ent{\Given{\Vect{Y}}{\Vect{X}}}}{\ent{\Vect{Y}}} \right)  \qquad \mbox{ and } \qquad \hS{\Vect{X}}{\Vect{Y}}  = \frac{\ent{\Given{\Vect{X}}{\Vect{Y}}} + \ent{\Given{\Vect{Y}}{\Vect{X}}}}{ \ent{\Vect{X}} + \ent{\Vect{Y}}}.
$$
which can be rewritten as \nibd[s:] 
\begin{eqnarray}
\hE{\Vect{X}}{\Vect{Y}} &=& 1-\frac{\IM{\Vect{X}}{\Vect{Y}}}{\pen[D]{\Vect{X}}{\Vect{Y}}}\quad \mbox{ with }\quad \pen[D]{\Vect{X}}{\Vect{Y}}=\left(\frac12\left(
\frac1{\ent{\Vect{X}}} + \frac1{\ent{\Vect{Y}}}\right) \right)^{-1}\label{def-hE}\\
\hS{\Vect{X}}{\Vect{Y}}  &=& 1-\frac{\IM{\Vect{X}}{\Vect{Y}}}{\pen[S]{\Vect{X}}{\Vect{Y}}} \quad \mbox{ with }\quad \pen[S]{\Vect{X}}{\Vect{Y}}=\frac12 \left(\ent{\Vect{X}} + \ent{\Vect{Y}} \right). \label{def-hS}
\end{eqnarray}
Their non normalized version are expressed as
$\HE{\Vect{X}}{\Vect{Y}} =  \pen[D]{\Vect{X}}{\Vect{Y}}-\IM{\Vect{X}}{\Vect{Y}}$ and 
$\HS{\Vect{X}}{\Vect{Y}} = \pen[S]{\Vect{X}}{\Vect{Y}}-\IM{\Vect{X}}{\Vect{Y}}$.

In this paper, we are interested in a large family of \ibd or \nibd with \compl[s] of the form:
\begin{equation} \label{eq-NIBg}
\pen[\alpha]{\Vect{X}}{\Vect{Y}}=g^{-1}\Big( \alpha \times g({m_{\Vect{X},\Vect{Y}}}) + (1-\alpha) \times g({M_{\Vect{X},\Vect{Y}}})  \Big)
\end{equation}
with ${m_{\Vect{X},\Vect{Y}}}={\min\left(\ent{\Vect{X}} , \ent{\Vect{Y}} \right)}$ and ${M_{\Vect{X},\Vect{Y}}}={\max\left(\ent{\Vect{X}} , \ent{\Vect{Y}} \right)}$ and where $0\leq\alpha < 1$ and $g(\cdot)$ is any monotone function on $\RR^+$. When it is not ambiguous we set $m={m_{\Vect{X},\Vect{Y}}}$ and $M={M_{\Vect{X},\Vect{Y}}}$. To be convinced that \ibd[s] and \nibd[s] with \compl[s] of the form~(\ref{eq-NIBg}) satisfy $(ii)$ of Definition~\ref{def-infoMes}, let us notice that
$$
\IM{\Vect{X}}{\Vect{Y}}= g^{-1}\left( \alpha g(\IM{\Vect{X}}{\Vect{Y}}) + (1-\alpha) g(\IM{\Vect{X}}{\Vect{Y}})  \right) \leq g^{-1}\left( \alpha g(m) + (1-\alpha) g(M)  \right).
$$
When $\alpha=0$, the \compl \penempty[\alpha] corresponds to \penempty[I]. When $\alpha=1$ the \compl defined by ${\min\left(\ent{\Vect{X}} , \ent{\Vect{Y}} \right)}$ and denoted by $\pen[\min]{\Vect{X}}{\Vect{Y}}$ does not satisfy  $(ii)$ of Definition~\ref{def-infoMes} and then \textbf{[P3]}. The associated  \IBempty[\min] (resp.  \NIBempty[\min]) is not an \ibd (resp. a \nibd[]). 

We pay now particular attention on the \compl[s] \penempty[D,\alpha], \penempty[S,\alpha], \penempty[R,\alpha] and \penempty[P,\alpha] of the form~(\ref{eq-NIBg}) respectively with $g^D(\cdot)=1/\cdot$, $g^S(\cdot)=\cdot$, $g^R(\cdot)=\sqrt{\cdot}$ and $g^P(\cdot)=\log(\cdot)$:
\begin{eqnarray}
\pen[D,\alpha]{\Vect{X}}{\Vect{Y}} &=& \left(\alpha \frac1{{\min\left(\ent{\Vect{X}} , \ent{\Vect{Y}} \right)}}+(1-\alpha)\frac1{{\max\left(\ent{\Vect{X}} , \ent{\Vect{Y}} \right)}}\right)^{-1}  \label{def-CD}\\
\pen[S,\alpha]{\Vect{X}}{\Vect{Y}} &=&  \alpha {\min\left(\ent{\Vect{X}} , \ent{\Vect{Y}} \right)} + (1-\alpha) {\max\left(\ent{\Vect{X}} , \ent{\Vect{Y}} \right)}. \label{def-CS}\\
\pen[R,\alpha]{\Vect{X}}{\Vect{Y}} &=& \left(\alpha\sqrt{{\min\left(\ent{\Vect{X}} , \ent{\Vect{Y}} \right)}}+(1-\alpha)\sqrt{{\max\left(\ent{\Vect{X}} , \ent{\Vect{Y}} \right)}}\right)^2  \label{def-CR}\\
\pen[P,\alpha]{\Vect{X}}{\Vect{Y}} &=&  {\min\left(\ent{\Vect{X}} , \ent{\Vect{Y}} \right)}^\alpha{\max\left(\ent{\Vect{X}} , \ent{\Vect{Y}} \right)}^{1-\alpha}. \label{def-CP}
\end{eqnarray}

The previous measures \IBempty[S],\NIBempty[S], \IBempty[D] and \NIBempty[D] are particular examples of such a family since the value of $\alpha=\frac 12$ leads to $\pen[1/2]{\Vect{X}}{\Vect{Y}}=g^{-1}\left(\frac12 g(\ent{\Vect{X}})+\frac12 g(\ent{\Vect{Y}}) \right)$.
When $\alpha=\frac12$,  \IBempty[\bullet,\alpha] and \NIBempty[\bullet,\alpha] will be simply denoted by \IBempty[\bullet] and \NIBempty[\bullet] where $\bullet$ stands for $S,R,P$ and $D$.

Let us first comment the particular expressions of the divergences \IBempty[S,\alpha] and \NIBempty[D,\alpha] associated to \penempty[D,\alpha] and  \penempty[S,\alpha] given by:
\begin{eqnarray*}
\IB[S,\alpha]{\Vect{X}}{\Vect{Y}}&=&\alpha \min \left( \ent{\Given{\Vect{X}}{\Vect{Y}}}, \ent{\Given{\Vect{Y}}{\Vect{X}}} \right) + (1-\alpha) \max \left(  \ent{\Given{\Vect{X}}{\Vect{Y}}}, \ent{\Given{\Vect{Y}}{\Vect{X}}} \right)\\
& =& \alpha \IB[\min]{\Vect{X}}{\Vect{Y}} + (1-\alpha) \DIB[I]{\Vect{X}}{\Vect{Y}}\\
\NIB[D,\alpha]{\Vect{X}}{\Vect{Y}}&=& \alpha \min \left( \frac{\ent{\Given{\Vect{X}}{\Vect{Y}}} }{\ent{\Vect{X}}}, \frac{\ent{\Given{\Vect{Y}}{\Vect{X}}} }{\ent{\Vect{Y}}} \right) + (1-\alpha) \max \left( 
\frac{\ent{\Given{\Vect{X}}{\Vect{Y}}} }{\ent{\Vect{X}}}, \frac{\ent{\Given{\Vect{Y}}{\Vect{X}}} }{\ent{\Vect{Y}}}\right) \\
&=& \alpha \NIB[\min]{\Vect{X}}{\Vect{Y}} + (1-\alpha) \dIB[I]{\Vect{X}}{\Vect{Y}}
\end{eqnarray*}

Clearly, the previous representation of \IB[S,\alpha]{\Vect{X}}{\Vect{Y}} (resp. \NIB[D,\alpha]{\Vect{X}}{\Vect{Y}}) as a convex combination of \IB[\min]{\Vect{X}}{\Vect{Y}} and  \DIB[I]{\Vect{X}}{\Vect{Y}} (resp. \NIB[\min]{\Vect{X}}{\Vect{Y}} and  \dIB[I]{\Vect{X}}{\Vect{Y}}) introduces a degree of freedom that could be useful for practical purposes in prediction framework where $\Vect{Y}$ could represent some  class variable. According to the parameter $\alpha$ one may favour to take into account between one or two prediction terms amoung $\ent{\Given{\Vect{X}}{\Vect{Y}}}$ and $\ent{\Given{\Vect{Y}}{\Vect{X}}}$ (resp. $\frac{\ent{\Given{\Vect{X}}{\Vect{Y}}} }{\ent{\Vect{X}}}$ and $\frac{\ent{\Given{\Vect{Y}}{\Vect{X}}} }{\ent{\Vect{Y}}}$). This possibility to introduce a non uniform mixing of the entropic contributions in the expression of the complexity terms seems to be not feasible by a direct adaptation of $\pen[I]{\Vect{X}}{\Vect{Y}}$.

\begin{remark}
By choosing $g(\cdot)=(\cdot)^\gamma$ for some $\gamma>0$, the \compl is given by $\pen[\gamma,\alpha]{\Vect{X}}{\Vect{Y}}= || \left(\alpha^{\frac1\gamma} m, (1-\alpha)^{\frac1\gamma} M\right) ||_{\gamma}$, where $||\Vect{x}||_{\gamma}=\left(\sum_{i=1}^2|x_i|^\gamma\right)^{1/\gamma}$ denotes the norm of some vector $\Vect{x}$ of length $2$. Note that for any $0\leq \alpha \leq 1$, we have
$$
({\alpha^\wedge})^{\frac1\gamma} ||(\ent{\Vect{X}},\ent{\Vect{Y}}) ||_{\gamma} \leq  \pen[\gamma,\alpha]{\Vect{X}}{\Vect{Y}} \leq 
({\alpha^{\vee}})^{\frac1\gamma} ||(\ent{\Vect{X}},\ent{\Vect{Y}}) ||_{\gamma} ,
$$
with ${\alpha^\wedge}=\min(\alpha,1-\alpha)$ and ${\alpha^{\vee}}=\max(\alpha,1-\alpha)$. When $\gamma$ goes to infinity $\pen[\gamma,\alpha]{\Vect{X}}{\Vect{Y}}$ converges towards $\pen[I]{\Vect{X}}{\Vect{Y}}$.
\end{remark}
 
\begin{remark}
The \compl $\penempty[\alpha]$ is invariant under linear transformation of $g$. In particular, $g$ and $-g$ provide the same \compl[.] Consequently, without loss of generality we could restrict $g$ to be an increasing function. 
\end{remark}

 Let us now propose a result to arrange these different examples considered in this paper. Before, some preliminary result is given.

\begin{lemma} \label{lem-C1C2}
Let $\penempty[(1)]$ and $\penempty[(2)]$ two \compl[s] of the form~(\ref{eq-NIBg}) with function $g_1$ and $g_2$. Assume either that the function $g_1 \circ g_2^{-1}$ is concave or that the function $g_2 \circ g_1^{-1}$ is convex, then $\pen[(1)]{\Vect{X}}{\Vect{Y}}\leq \pen[(2)]{\Vect{X}}{\Vect{Y}}$
\end{lemma}

\begin{proof}
By rewriting $g_1=(g_1\circ g_2^{-1})\circ g_2 $ when $g_1\circ g_2^{-1}$ is concave and $g_1^{-1}=g_2^{-1}\circ (g_2\circ g_1^{-1})$ when $(g_2\circ g_1^{-1})$ is convex, one may assert
\begin{eqnarray}
g_1^{-1}(\alpha g_1(m)+ (1-\alpha)g_1(M) ) &\leq& \left\{
\begin{array}{l}
g_2^{-1}(\alpha (g_2\circ g_1^{-1})\circ g_1(m) + (1-\alpha) (g_2\circ g_1^{-1})\circ g_1(M)) \\
g_1^{-1}( g_1\circ g_2^{-1} (\alpha g_2(m)+ (1-\alpha)g_2(M) ))
\end{array}
\right. \nonumber \\
&\leq & g_2^{-1}(\alpha g_2(m)+ (1-\alpha)g_2(M) )  \nonumber
\end{eqnarray}
where $m={\min\left(\ent{\Vect{X}} , \ent{\Vect{Y}} \right)}$ and $M={\max\left(\ent{\Vect{X}} , \ent{\Vect{Y}} \right)}$.
\end{proof}

\begin{proposition}
For any $\IBempty[(1)]$, $\IBempty[(2)]$ \ibd[s] or any $\NIBempty[(1)]$, $\NIBempty[(2)]$ \nibd[s] with respective \compl[s] $\penkempty{1}$ and $\penkempty{2}$, the following equivalence holds:
\begin{equation} \label{relationOrdre}
\IBk{1}{\Vect{X}}{\Vect{Y}} \leq \IBk{2}{\Vect{X}}{\Vect{Y}} \Longleftrightarrow \NIBk{1}{\Vect{X}}{\Vect{Y}} \leq \NIBk{2}{\Vect{X}}{\Vect{Y}} \Longleftrightarrow  \penk{1}{\Vect{X}}{\Vect{Y}}  \leq \penk{2}{\Vect{X}}{\Vect{Y}}.
\end{equation}
Since, for any $0\leq\alpha\leq \alpha^\prime\leq 1$,
\begin{equation} \label{eq-CalphaCI}
\pen[\alpha]{\Vect{X}}{\Vect{Y}} \leq \pen[I]{\Vect{X}}{\Vect{Y}}\quad\mbox{ and }\quad \pen[\alpha]{\Vect{X}}{\Vect{Y}} \leq \pen[\alpha^\prime]{\Vect{X}}{\Vect{Y}}
\end{equation}
the associated \ibd[s] and \nibd[s] are then ordered according to equation~(\ref{relationOrdre}).
Furthermore, a similar result holds for the main examples of this paper since
\begin{equation} \label{inegDist}
\pen[D,\alpha]{\Vect{X}}{\Vect{Y}} \leq \pen[P,\alpha]{\Vect{X}}{\Vect{Y}} \leq \pen[R,\alpha]{\Vect{X}}{\Vect{Y}} \leq \pen[S,\alpha]{\Vect{X}}{\Vect{Y}} \leq \pen[I]{\Vect{X}}{\Vect{Y}} \leq \pen[E]{\Vect{X}}{\Vect{Y}} 
\end{equation}

\end{proposition}

\begin{proof}
Equation (\ref{relationOrdre}) is direct. The left-hand side of (\ref{eq-CalphaCI}) comes from
$$
\pen[\alpha]{\Vect{X}}{\Vect{Y}} = g^{-1} \left(\alpha g({\min\left(\ent{\Vect{X}} , \ent{\Vect{Y}} \right)})+ (1-\alpha)g({\max\left(\ent{\Vect{X}} , \ent{\Vect{Y}} \right)})\right) \leq g^{-1} \left( g \left( {\max\left(\ent{\Vect{X}} , \ent{\Vect{Y}} \right)} \right) \right) = \pen[I]{\Vect{X}}{\Vect{Y}},
$$
and the right-hand side is direct. Since $g^P \circ (g^D)^{-1}(\cdot)= - \log(\cdot)$, $g^R \circ (g^P)^{-1}(\cdot)= \exp(\frac12 \cdot)$ and $g^S \circ (g^R)^{-1}(\cdot)= (\cdot)^2$ are convex functions, (\ref{inegDist}) is a direct consequence of Lemma~\ref{lem-C1C2}.
\end{proof}

\begin{remark}
By assuming either that $g(\cdot)$ is a convex function or that $g^{-1}(\cdot)$ is a concave function, the following inequality holds
$$
\pen[\alpha]{\Vect{X}}{\Vect{Y}} \leq \alpha m + (1- \alpha) M = \pen[S,\alpha]{\Vect{X}}{\Vect{Y}}
$$
which means that any \IBempty[\alpha] (resp. \NIBempty[\alpha]) (satisfying the previous assumption) is upper bounded by \IBempty[S,\alpha] (resp. \NIBempty[S,\alpha]). 
\end{remark}


The following proposition gives a larger class of examples of \ibd[s] and \nibd[s.]

\begin{proposition} \label{prop-ExMesInfo}
Let $(\alpha^{(j)})_{j=1,\ldots,J}$ be some vector of probability weights for some $J\geq 1$.

$(i)$ Let $\NIBkempty{1}, \ldots, \NIBkempty{J}$, $J$ \nibd[s,] then the measure defined by
\begin{equation} \label{ex-CCMesInfo}
\NIB{\Vect{X}}{\Vect{Y}} = \sum_{j=1}^J \alpha^{(j)} \;\NIBk{j}{\Vect{X}}{\Vect{Y}}
\end{equation}
is a \nibd with \compl given by 
\begin{equation}\label{eq-normCCDist}
\pen{\Vect{X}}{\Vect{Y}} = \left( \sum_{j=1}^J \frac{\alpha^{(j)}}{\penk{j}{\Vect{X}}{\Vect{Y}}}\right)^{-1}.
\end{equation}

$(ii)$ Let $\IBkempty{1}, \ldots, \IBkempty{j}$, $J$ \ibd[s] and $\NIBkempty{1}, \ldots, \NIBkempty{j}$, $J$ \nibd[s] with \compl[s] $\penk{1}{\Vect{X}}{\Vect{Y}}, \ldots, \penk{J}{\Vect{X}}{\Vect{Y}}$ then the measures defined by
\begin{equation} \label{ex-CCNorm}
\IB{\Vect{X}}{\Vect{Y}} = \pen{\Vect{X}}{\Vect{Y}} -\IM{\Vect{X}}{\Vect{Y}} \quad \mbox{ and } \quad 
\NIB{\Vect{X}}{\Vect{Y}} = 1 -\frac{\IM{\Vect{X}}{\Vect{Y}}}{\pen{\Vect{X}}{\Vect{Y}}}, \qquad \mbox{ with } \qquad \pen{\Vect{X}}{\Vect{Y}}= \sum_{j=1}^J \alpha^{(j)} \; \penk{j}{\Vect{X}}{\Vect{Y}}
\end{equation}
are also respectively an \ibd and a \nibd[.]
\end{proposition}

The proof is immediate.

\subsection{Around the property \textbf{[P3]}}

The fact that an \ibd \IBempty (resp. \nibd \NIBempty) satisfies the property \textbf{[P3]} may be expressed by: $\IB{\Vect{X}}{\Vect{Y}}=0 \Leftrightarrow\DIB[I]{\Vect{X}}{\Vect{Y}}=0$ (resp. $\NIB{\Vect{X}}{\Vect{Y}}=0 \Leftrightarrow\dIB[I]{\Vect{X}}{\Vect{Y}}=0$). In fact, \textbf{[P3]}
should be extended to the more useful assumption: $\IB{\Vect{X}}{\Vect{Y}}$ (or $\NIB{\Vect{X}}{\Vect{Y}}$) is near from  minimum 0 if and only if $\Vect{X}$ and $\Vect{Y}$ share almost the same information. This may be translated by the following implications related to an \ibd \IBempty (resp. a \nibd \NIBempty):
\begin{itemize}
\item for all $\gamma>0$ there exists $\varepsilon>0$ such that for all $(\Vect{X}, \Vect{Y})\in\Upsilon$
$$ 
\IB{\Vect{X}}{\Vect{Y}}\leq \varepsilon \Longrightarrow  \DIB[I]{\Vect{X}}{\Vect{Y}}\leq \gamma \quad (\mbox{resp. } \quad
\NIB{\Vect{X}}{\Vect{Y}}\leq \varepsilon \Longrightarrow  \dIB[I]{\Vect{X}}{\Vect{Y}}\leq \gamma).
$$
\item for all  $\varepsilon>0$ there exists  $\gamma>0$ such that for all $(\Vect{X}, \Vect{Y})\in \Upsilon$
$$
\DIB[I]{\Vect{X}}{\Vect{Y}}\leq \gamma \Longrightarrow  \IB{\Vect{X}}{\Vect{Y}}\leq \varepsilon
\quad (\mbox{resp. } \quad
\dIB[I]{\Vect{X}}{\Vect{Y}}\leq \gamma \Longrightarrow  \NIB{\Vect{X}}{\Vect{Y}}\leq \varepsilon).
$$
\end{itemize}
An \ibd \IBempty (resp. a \nibd \NIBempty) inherits of the previous property if it satisfies:

\textbf{[P3bis($\Upsilon,k_1,k_2$)]} there exists some positive constants $k_1, k_2$ ($k_1\leq k_2$) such that for all $\sloppy{ (\Vect{X},\Vect{Y})\in \Upsilon\subset\Gamma^2}$:
\begin{equation} \label{eq-P3bis}
k_1\; \DIB[I]{\Vect{X}}{\Vect{Y}} \leq \IB{\Vect{X}}{\Vect{Y}} \leq k_2 \; \DIB[I]{\Vect{X}}{\Vect{Y}}\quad (\mbox{resp. } \quad k_1\; \dIB[I]{\Vect{X}}{\Vect{Y}} \leq \NIB{\Vect{X}}{\Vect{Y}} \leq k_2 \; \dIB[I]{\Vect{X}}{\Vect{Y}}).
\end{equation}

Among our examples, we assert that \DEempty and  \dEempty  both satisfy \textbf{[P3bis($\Gamma^2,1,2$)]} that is
$$ 
\DIB[I]{\Vect{X}}{\Vect{Y}} \leq \DIB[E]{\Vect{X}}{\Vect{Y}} \leq 2 \DIB[I]{\Vect{X}}{\Vect{Y}} \quad (\mbox{resp. } \dIB[I]{\Vect{X}}{\Vect{Y}} \leq \dIB[E]{\Vect{X}}{\Vect{Y}} \leq 2 \dIB[I]{\Vect{X}}{\Vect{Y}}).
$$

Most of \compl[s] considered in this paper are of the particular form~(\ref{eq-NIBg}) where the function $g(\cdot)$ is a monotone function on $\RR^+$. 
From~(\ref{eq-CalphaCI}), we can point out that for such \compl[s] (expressed in terms of \IBempty or \NIBempty), the constant $k_2$ is equal to~1. Moreover, we assert that if \IBempty satisfies \textbf{[P3bis($\Upsilon,k_1,1$)]} then the associated \NIBempty  also satisfies \textbf{[P3bis($\Upsilon,k_1,1$)]} since 
$$
k_1  \;\dIB[I]{\Vect{X}}{\Vect{Y}} = \frac{k_1\DIB[I]{\Vect{X}}{\Vect{Y}}}{\pen[I]{\Vect{X}}{\Vect{Y}}} \leq \frac{\IB{\Vect{X}}{\Vect{Y}}}{\pen{\Vect{X}}{\Vect{Y}}} =\NIB{\Vect{X}}{\Vect{Y}}.
$$
And so in the rest of this section, the results presented hereafter will be only expressed for \ibd[s.]

Furthermore, we now consider only \compl[s] of the form~(\ref{eq-NIBg}) defined through a function $g(\cdot)$ continuously differentiable on some set $\mathcal{D}^g\subset \RR^+$. Let us first introduce the two following subsets of $\mathcal{D}^g$:
$$
\mathcal{E}_1^g = \left\{ \Theta\subset \mathcal{D}^g: 0< \kappa_{\inf,\Theta}^{g}<\kappa_{\sup,\Theta}^{g}<+\infty\right\}
\quad \mbox{ and } \quad 
\mathcal{E}_2^{g,\alpha}= \left\{ \Theta\subset \mathcal{E}_1^g: \frac{\alpha \kappa_{\sup,\Theta}^{g}}{  \kappa_{\inf,\Theta}^{g}}<1\right\},
$$
with $\kappa_{\inf,\Theta}^{g}= \inf_{x \in \Theta} |g^\prime(x)|$ and $\kappa_{\sup,\Theta}^{g} = \sup_{x \in \Theta} |g^\prime(x)|$. Denote also by ${\alpha^\wedge}=\min(\alpha,1-\alpha)$.

In the sequel, two  results ensuring that an \ibd \IBempty[\alpha] of the 
form~(\ref{eq-NIBg}) satisfies \textbf{[P3bis($\Upsilon,k_1,k_2$)]}, 
are proposed. The difference relies upon the framework: the constants $k_1$ and $k_2$ differ whenever the set $\Upsilon$ differs.

\begin{proposition} \label{prop1-P3bis} 
For any $\Theta \in  \mathcal{E}_1^g $ the \ibd \IBempty[\alpha] satisfies \textbf{[P3bis($\Upsilon_\Theta, {\alpha^\wedge}\frac{\kappa_{\inf,\Theta}^{g}}{\kappa_{\sup,\Theta}^{g}},1$)]}  with \sloppy{$\Upsilon_\Theta= \left\{ (\Vect{X},\Vect{Y}) \in \Gamma^2: \ent{\Vect{X}},\ent{\Vect{Y}},\IM{\Vect{X}}{\Vect{Y}} \in \Theta \right\}.$}
\end{proposition}

\begin{proof}
Denote by $x=\min \left(\ent{\Given{\Vect{X}}{\Vect{Y}}}, \ent{\Given{\Vect{Y}}{\Vect{X}}} \right)$, $y=  \max \left(\ent{\Given{\Vect{X}}{\Vect{Y}}}, \ent{\Given{\Vect{Y}}{\Vect{X}}} \right)  $ and $z=\IM{\Vect{X}}{\Vect{Y}}$. There exists $c_1,c_2,c_3$ such that 
\begin{eqnarray*}
g^{-1} \left( \frac{g(x+z) +g(y+z)}2 \right) - z &=& 
 \left( \alpha (g(x+z)-g(z)) + (1-\alpha) (g(y+z)-g(z))\right) (g^{-1})^\prime(c_1) \\
&=& \alpha |g^\prime(c_2)| |(g^{-1})^\prime(c_1)| \times x + (1-\alpha) |g^\prime(c_3)| |(g^{-1})^\prime(c_1)| \times y,
\end{eqnarray*}
with $c_1\in\left[\min\left(g(z),\alpha g(x+z)+(1-\alpha)g(y+z)\right), \max\left(g(z),\alpha g(x+z)+(1-\alpha)g(y+z)\right) \right]$, $c_2\in[z,x+z]$ and $c_3\in[z,y+z]$. Then, we obtain for all $x,y,z$:
$$
g^{-1} \left( \frac{g(x+z) +g(y+z)}2 \right) - z \geq 
{\alpha^\wedge} \frac{\kappa_{\inf,\Theta}^{g}}{ \kappa_{\sup,\Theta}^{g}} \max(x,y) 
$$
which means that $ {\alpha^\wedge} \frac{\kappa_{\inf,\Theta}^{g}}{ \kappa_{\sup,\Theta}^{g}}  \;\DIB[I]{\Vect{X}}{\Vect{Y}} \leq \IB[\alpha]{\Vect{X}}{\Vect{Y}}$.
\end{proof}

\begin{proposition} \label{prop2-P3bis}
For any $\Theta \in  \mathcal{E}_2^g $ the \ibd \IBempty[\alpha] satisfies \textbf{[P3bis($\Gamma^2_\Theta, 1-\alpha\frac{\kappa_{\sup,\Theta}^{g}}{\kappa_{\inf,\Theta}^{g}},1$)]}  with 
$\Gamma_\Theta=\left\{\Vect{Z}\in \Gamma: \ent{\Vect{Z}} \in \Theta \right\}$. 
\end{proposition}

\begin{proof}
\begin{eqnarray*}
\DIB[I]{\Vect{X}}{\Vect{Y}}-\IB[\alpha]{\Vect{X}}{\Vect{Y}} = \pen[I]{\Vect{X}}{\Vect{Y}}-\pen{\Vect{X}}{\Vect{Y}} &=& \alpha (g^{-1})^\prime(c_1) \left( g(\max(\ent{\Vect{X}},\ent{\Vect{Y}}))- g(\min(\ent{\Vect{X}},\ent{\Vect{Y}}))\right) \\
&=& \alpha|(g^{-1})^\prime(c_1)| |g^\prime(c_2)| \left| \ent{\Vect{X}}-\ent{\Vect{Y}}\right|, 
\end{eqnarray*}
with $c_1\in \left[ g(\min(\ent{\Vect{X}},\ent{\Vect{Y}})),g(\max(\ent{\Vect{X}},\ent{\Vect{Y}}))\right]$ and $c_2 \in \left[ \min(\ent{\Vect{X}},\ent{\Vect{Y}}), \max(\ent{\Vect{X}},\ent{\Vect{Y}})\right]$. Then we obtain   
$$
\DIB[I]{\Vect{X}}{\Vect{Y}}-\IB[\alpha]{\Vect{X}}{\Vect{Y}} \leq \alpha \frac{\kappa_{\sup,\Theta}^{g}}{\kappa_{\inf,\Theta}^{g}} \times \DIB[I]{\Vect{X}}{\Vect{Y}}
$$
which leads to the result.
\end{proof}

For sake of simplicity, we denote by $\kappa_{\inf,\Theta}^\bullet$ and $\kappa_{\sup,\Theta}^\bullet$ instead of $\kappa_{\inf,\Theta}^{g^\bullet}$ and $\kappa_{\sup,\Theta}^{g^\bullet}$

The following result is devoted to our different examples. We apply the 
two previous propositions and present a new result obtained by taking 
into account the specific form of each example.

\begin{proposition} \label{prop3-P3bis}
\IBempty[\bullet,\alpha] satisfies \textbf{[P3bis($\Upsilon_\Theta,k_1^{a,\bullet},1$)]} (from Proposition~\ref{prop1-P3bis}), \textbf{[P3bis($\Gamma^2_\Theta, k_{1}^{b,\bullet},1$)]} (from Proposition~\ref{prop2-P3bis}) and \textbf{[P3bis($\Gamma^2_\Theta, k_{1}^{c,\bullet},1$)]}
where $\bullet$ stands for $S,R,P$ and $D$, and
\begin{center}
\begin{tabular}{|c|c|c|c|c|c|c|}
\hline 
$\bullet$& $\Theta$ &$ \kappa_{\inf,\Theta}^\bullet  $ & $ \kappa_{\sup,\Theta}^\bullet\quad $ &$ k_1^{a,\bullet}={\alpha^\wedge}\frac{\kappa_{\inf,\Theta}}{\kappa_{\sup,\Theta}^\bullet} $
&$k_1^{b,\bullet}=1-\alpha\frac{\kappa_{\sup,\Theta}}{\kappa_{\inf,\Theta}^\bullet} $ &
$k_1^{c,\bullet}$
\\ \hline
$S$ & $\RR^+$ & 1 & 1 & ${\alpha^\wedge}$&$1-\alpha$ &  \\ \hline 
$R$ & $[c_1,c_2]$ &$\frac1{2\sqrt{c_2}}$ & $\frac1{2\sqrt{c_1}}$ &
$\frac{{\alpha^\wedge}}{\sqrt{\rho}}$& $1-\alpha\sqrt{\rho}$ (if $\rho<\frac1{\alpha^2}$)&$(1-\alpha)\left(1-\frac\alpha{(1+\frac1{\sqrt{\rho}})^2}\right)$\\ \hline
$R$ & $\RR^+$ & & & & &$(1-\alpha)^2$\\ \hline
$P$& $[c_1,c_2]$ &$\frac1{c_2}$ & $\frac1{c_1}$ &
$\frac{{\alpha^\wedge}}{ \rho}$ & $1-\alpha\rho$ (if $\rho<\frac1\alpha $)&
$\frac{\rho^\alpha-1}{\rho-1}$
\\ \hline
$D$ &$[c_1,c_2]$ & $\frac1{c_2^2}$ & $\frac1{c_1^2}$ &
$\frac{{\alpha^\wedge}}{ \rho^2}$ & $1-\alpha\rho^2$ (if $\rho<\frac1{\sqrt{\alpha}}$)&
$\frac1{1+\frac\alpha{1-\alpha}\rho}$\\\hline
\end{tabular}
\end{center}
with $0< c_1\leq c_2<+\infty$, $\rho=\frac{c_2}{c_1}$.
\end{proposition}

\begin{proof}
The computations of $k_1^{a,\bullet}$ and $k_1^{b,\bullet}$ derive from Proposition~\ref{prop1-P3bis} and \ref{prop2-P3bis}. Hence, let us concentrate only on $k_1^{c,\bullet}$ for the complexity terms \penempty[R], \penempty[P] and \penempty[D]. Let us denote by $m={\min\left(\ent{\Vect{X}} , \ent{\Vect{Y}} \right)}$ and by $M={\max\left(\ent{\Vect{X}} , \ent{\Vect{Y}} \right)}$.
\begin{itemize}
\item \Compl $\penempty[R]$:
\begin{eqnarray*}
\DIB[I]{\Vect{X}}{\Vect{Y}} - \IB[R,\alpha]{\Vect{X}}{\Vect{Y}} &=& \alpha(1-\alpha) \left( \sqrt{M}-\sqrt{m} \right)^2 + \alpha(M-m) \\
&=& \alpha(1-\alpha) \frac{(M-m)^2}{(\sqrt{M}+\sqrt{m} )^2} + \alpha(M-m) \\
&\leq &  \alpha(1-\alpha)\frac{(\DIB[I]{\Vect{X}}{\Vect{Y}})^2}{(\sqrt{M}+\sqrt{m} )^2}  +\alpha\DIB[I]{\Vect{X}}{\Vect{Y}}\\
\end{eqnarray*}
And so, 
$$
\IB[R,\alpha]{\Vect{X}}{\Vect{Y}} \geq (1-\alpha) \DIB[I]{\Vect{X}}{\Vect{Y}} \left( 1- \alpha \frac{\DIB[I]{\Vect{X}}{\Vect{Y}}}{ (\sqrt{M}+\sqrt{m} )^2 } \right)
$$
The result is obtained by noticing that
\begin{eqnarray*}
\frac{\DIB[I]{\Vect{X}}{\Vect{Y}}}{\left( \sqrt{m} + \sqrt{M} \right)^2} &\leq&
\frac{M}{\left( \sqrt{m} + \sqrt{M} \right)^2 } = \frac1{\left(1+ \sqrt{\frac{m}{M}}\right)^2}
\\
&\leq& \frac{1}{\left(1+ \sqrt{\frac{c_1}{c_2}}\right)^2} \leq 1.
\end{eqnarray*}


\item \Compl \penempty[P]: by using a Taylor expansion with integral rest, one obtains

\begin{eqnarray*}
\DIB[I]{\Vect{X}}{\Vect{Y}} - \IB[P,\alpha]{\Vect{X}}{\Vect{Y}} &=& M^\alpha \left( M^{1-\alpha} - m^{1-\alpha}\right) \\
&=& M^\alpha (M-m) \times \int_0^1 \frac{1-\alpha}{\left(m+t(M-m)\right)^\alpha} dt\\
&\leq & (M-m) \int_0^1 \frac{1-\alpha}{\left(\frac1\rho + t(1-\frac1\rho) \right)^\alpha}dt \\
&\leq & \DIB[I]{\Vect{X}}{\Vect{Y}} \frac1{1-\frac1\rho} \left[ \left( \frac1\rho+t(1-\frac1\rho)\right)^{1-\alpha}\right]_0^1= \DIB[I]{\Vect{X}}{\Vect{Y}} \frac{1-\left(\frac1\rho\right)^{1-\alpha}}{1-\frac1\rho}
\end{eqnarray*}
And so, 
$$
\IB[P,\alpha]{\Vect{X}}{\Vect{Y}} \geq \DIB[I]{\Vect{X}}{\Vect{Y}} \left(1-  \frac{1-\left(\frac1\rho\right)^{1-\alpha}}{1-\frac1\rho}\right),
$$
which leads to the result.

\item \Compl \penempty[D]:
$$
\DIB[I]{\Vect{X}}{\Vect{Y}} - \IB[D,\alpha]{\Vect{X}}{\Vect{Y}} = M - \frac{m M}{\alpha M + (1-\alpha)m} = \frac{\alpha M}{\alpha M + (1-\alpha)m} (M-m) \leq \frac1{1+\frac{1-\alpha}{\alpha} \frac{c_1}{c_2}} \times \DIB[I]{\Vect{X}}{\Vect{Y}}
$$\end{itemize}
\end{proof}


\subsection{Around the triangular inequality's property} \label{sec-infoMesDist}

The question arises now whether an \ibd or a \nibd satisfies the property \textbf{[P6]} that is a triangular inequality. The following proposition establishes sufficient conditions for such measures to constitute a metric.

\begin{lemma} \label{lem-toolentJointIM}
\begin{eqnarray}
\ent{\Vect{X},\Vect{Y}} &\leq& \ent{ \Vect{X}, \Vect{Z}} +  \ent{ \Vect{Y}, \Vect{Z}} -\ent{\Vect{Z}}\label{eq-resEntJoint} \\
\IM{\Vect{X}}{\Vect{Y}} &\geq& \IM{\Vect{X}}{\Vect{Z}} +\IM{\Vect{Y}}{\Vect{Z}}-\ent{\Vect{Z}}\label{eq-resIMTri}
\end{eqnarray}
\end{lemma}
\begin{proof}
From general properties on entropy, one can obtain
\begin{equation}
\ent{\Vect{X},\Vect{Y}} \leq \ent{\Vect{X},\Vect{Y}, \Vect{Z}}= \ent{ \Vect{X}, \Vect{Z}} + \ent{ \Given{\Vect{Y}}{\Vect{X}, \Vect{Z}} }\leq \ent{ \Vect{X}, \Vect{Z}}+ \ent{ \Given{\Vect{Y}}{ \Vect{Z}} } = \ent{ \Vect{X}, \Vect{Z}} +  \ent{ \Vect{Y}, \Vect{Z}} -\ent{\Vect{Z}}.
\end{equation}
Equation~(\ref{eq-resIMTri}) directly derives from (\ref{eq-defIM}).
\end{proof}

\begin{proposition} \label{prop-inegTriInfoMes}
Assume the \compl defining an \ibd satisfies the following property: \\
\noindent $(i)$
\begin{equation} \label{hyp1inegTri}
\pen{\Vect{X}}{\Vect{Y}} \leq \pen{\Vect{X}}{\Vect{Z}} + \pen{\Vect{Y}}{\Vect{Z}}  - \ent{\Vect{Z}}.
\end{equation}
Then, the associated \ibd satisfies the triangular inequality, that is
\begin{equation} \label{eq-inegTriIB}
\IB{\Vect{X}}{\Vect{Y}} \leq \IB{\Vect{X}}{\Vect{Z}} + \IB{\Vect{Y}}{\Vect{Z}}.
\end{equation}
In addition, if $\penempty$ satisfies \\
$(ii)$
\begin{equation} \label{hyp2inegTri}
\pen{\Vect{X}}{\Vect{Z}}\geq \max\left( \ent{\Vect{X}},\ent{\Vect{Z}}\right),
\end{equation}
then the associated \nibd satisfies also a triangular inequality, that is
\begin{equation} \label{eq-inegTriInfoMes}
\NIB{\Vect{X}}{\Vect{Y}} \leq \NIB{\Vect{X}}{\Vect{Z}} + \NIB{\Vect{Y}}{\Vect{Z}}.
\end{equation}
\end{proposition}

\begin{proof}
Since the following quantity
$$
A = -(\pen{\Vect{X}}{\Vect{Y}} -\IM{\Vect{X}}{\Vect{Y}})  +
(\pen{\Vect{X}}{\Vect{Z}} - \IM{\Vect{X}}{\Vect{Z}}) +
(\pen{\Vect{Y}}{\Vect{Z}} - \IM{\Vect{Y}}{\Vect{Z}}),
$$
is nonnegative from (\ref{eq-resIMTri}) and (\ref{hyp1inegTri}), we have immediately~(\ref{eq-inegTriIB}). Moreover, the following equation is valid
\begin{equation} \label{eq-prInegTri1}
\NIB{\Vect{X}}{\Vect{Y}} \leq 1 - \frac{\IM{\Vect{X}}{\Vect{Y}}}{
\pen{\Vect{X}}{\Vect{Y}} + A},
\end{equation}
Now, it is also easy to see from~(\ref{hyp2inegTri}) that
$$
A + \pen{\Vect{X}}{\Vect{Y}} \geq \pen{\Vect{X}}{\Vect{Z}}+\pen{\Vect{Y}}{\Vect{Z}}- \ent{\Vect{Z}} \geq \max \left( \pen{\Vect{X}}{\Vect{Z}},\pen{\Vect{Y}}{\Vect{Z}}\right).
$$
From (\ref{eq-prInegTri1}) it follows
$$
\NIB{\Vect{X}}{\Vect{Y}} \leq \frac{ \pen{\Vect{X}}{\Vect{Z}} - \IM{\Vect{X}}{\Vect{Z}} +
\pen{\Vect{Y}}{\Vect{Z}} - \IM{\Vect{Y}}{\Vect{Z}} }{
\max \left( \pen{\Vect{X}}{\Vect{Z}},\pen{\Vect{Y}}{\Vect{Z}}\right)} \leq
\frac{\pen{\Vect{X}}{\Vect{Z}} - \IM{\Vect{X}}{\Vect{Z}}}{\pen{\Vect{X}}{\Vect{Z}}} +
\frac{\pen{\Vect{Y}}{\Vect{Z}} - \IM{\Vect{Y}}{\Vect{Z}}}{\pen{\Vect{Y}}{\Vect{Z}}}
= \NIB{\Vect{X}}{\Vect{Z}} + \NIB{\Vect{Y}}{\Vect{Z}}.
$$
\end{proof}

\begin{remark} \label{rem-metric1}
In Proposition~\ref{prop-inegTriInfoMes}, there is no implication between~(\ref{hyp1inegTri}) and~(\ref{hyp2inegTri}). Indeed, one may check that the \nibd $\NIBempty[S]$ (with $\alpha=1/2$ for example)  satisfies the first one but not the second one. Now consider a \nibd with \compl $\pen{\Vect{X}}{\Vect{Y}}=\max\left(\ent{\Vect{X}},\ent{\Vect{Y}}\right)+ \ent{\Given{\Vect{X}}{\Vect{Y}}} \ent{\Given{\Vect{Y}}{\Vect{X}}}$. By choosing $\Vect{X}, \Vect{Y}$ and $\Vect{Z}$ such that $ \ent{\Given{\Vect{X}}{\Vect{Y}}}= \ent{\Given{\Vect{Y}}{\Vect{X}}}=\IM{\Vect{X}}{\Vect{Y}}=\ent{\Vect{Z}}/3=\ent{\Vect{X},\Vect{Y}}/3>2$, one asserts that~(\ref{hyp2inegTri}) is satisfied but not~(\ref{hyp1inegTri}).
\end{remark}

\begin{remark} \label{rem-metric2}
Let us consider a \nibd $\NIBempty$ with \compl
$\pen{\Vect{X}}{\Vect{Y}} = \pen[\prime]{\Vect{X}}{\Vect{Y}} + \max \left( \ent{\Vect{X}}, \ent{\Vect{Y}}\right)$ such that $\pen[\prime]{\Vect{X}}{\Vect{Y}}\geq 0$ (necessarily $\pen[\prime]{\Vect{X}}{\Vect{Y}}= 0$ whenever $\Vect{X} \sim\Vect{Y}$). Then, $\IBempty$ and $\NIBempty$ satisfy a triangular inequality if $\penempty[\prime]$  also satisfies  a triangular inequality. However, this is not a necessary condition. Indeed, the triangular inequality is not satisfied for the same example of the previous remark with $\pen[\prime]{\Vect{X}}{\Vect{Y}}=\ent{\Given{\Vect{X}}{\Vect{Y}}} \ent{\Given{\Vect{Y}}{\Vect{X}}}$ for which $\pen[\prime]{\Vect{X}}{\Vect{Z}}=\pen[\prime]{\Vect{Y}}{\Vect{Z}}=0$ whereas $\pen[\prime]{\Vect{X}}{\Vect{Y}} > 0$.
\end{remark}

Let us now propose some examples and consequences through the following corollary.

\begin{corollary} \label{cor-MesInfoDist} ${ }$

$(i)$ The measures $\DEempty$, $\DIempty$ satisfy the condition~(\ref{hyp1inegTri}) and so are metrics.

$(ii)$ The measures $\dEempty$ and $\dIempty$ satisfy the conditions~(\ref{hyp1inegTri}) and~(\ref{hyp2inegTri}) and so are metrics.

$(iii)$ The measure $\DIBempty[S,\alpha]$ for $\alpha\leq \frac12$  satisfies the condition~(\ref{hyp1inegTri}) and so is a metric. Moreover, when $\alpha>\frac12$, this measure does not satisfy~(\ref{hyp1inegTri}).

$(iv)$ Let $(\alpha^{(j)})_{j=1,\ldots,J}$ be some vector of probability weights for some $J\geq 1$. Let $\IBkempty{1}, \ldots, \IBkempty{J}$, $J$ \ibd[s] (resp. $\NIBkempty{1}, \ldots, \NIBkempty{J}$, $J$ \nibd[s]) with \compl[s] $\penk{1}{\Vect{X}}{\Vect{Y}}, \ldots, \penk{J}{\Vect{X}}{\Vect{Y}}$ satisfying~(\ref{hyp1inegTri}) (resp. (\ref{hyp1inegTri}) and~(\ref{hyp2inegTri})) then these measures defined by~(\ref{ex-CCNorm}) satisfy a triangular inequality.

\end{corollary}

\begin{proof}
$(i)$ and $(ii)$ Equation (\ref{eq-resEntJoint}) corresponds exactly to (\ref{hyp1inegTri}) for $\pen[E]{\Vect{X}}{\Vect{Y}}=\ent{\Vect{X},\Vect{Y}}$. And since $\ent{\Vect{X},\Vect{Z}}\geq \max\left( \ent{\Vect{X}}, \ent{\Vect{Z}}\right)$, we have proved that $\DEempty$ and $\dEempty$ are metrics. Concerning $\DIempty$ and $\dIempty$, the \compl corresponds to $\pen[I]{\Vect{X}}{\Vect{Y}}=\max\left( \ent{\Vect{X}},\ent{\Vect{Y}}\right)$. Thus it is sufficient to prove~(\ref{hyp1inegTri}) which is quite obvious. Indeed,
$$
\max\left( \ent{\Vect{X}}, \ent{\Vect{Z}}\right) + \max\left( \ent{\Vect{Y}}, \ent{\Vect{Z}}\right) - \ent{\Vect{Z}} \geq \max\left( \ent{\Vect{X}}, \ent{\Vect{Y}}\right).
$$

$(iii)$ Let $m={\min\left(\ent{\Vect{X}} , \ent{\Vect{Y}} \right)}$ and $M={\max\left(\ent{\Vect{X}} , \ent{\Vect{Y}} \right)}$. We distinguish three cases~:
\begin{itemize}
\item $\ent{\Vect{Z}}< m$:
$$
\pen[S,\alpha]{\Vect{X}}{\Vect{Z}}+\pen[S,\alpha]{\Vect{Y}}{\Vect{Z}}-\ent{\Vect{Z}}=(2\alpha-1)\ent{\Vect{Z}}+(1-\alpha)(m+M)
$$
If $\alpha>\frac12$ and $\ent{\Vect{X}}=\ent{\Vect{Y}}$, the right-hand side of the previous equation equals $(1-2\alpha)(\pen[S,\alpha]{\Vect{X}}{\Vect{Y}}-\ent{\Vect{Z}})+\pen[S,\alpha]{\Vect{X}}{\Vect{Y}} < \pen[S,\alpha]{\Vect{X}}{\Vect{Y}} $. And so, (\ref{hyp1inegTri}) can never be satisfied for $\alpha>\frac12$. Now, if $\alpha\leq \frac12$, we have
$$
\pen[S,\alpha]{\Vect{X}}{\Vect{Z}}+\pen[S,\alpha]{\Vect{Y}}{\Vect{Z}}-\ent{\Vect{Z}}>(1-\alpha)(m+M)\geq \pen[S,\alpha]{\Vect{X}}{\Vect{Y}}
$$
\item $\ent{\Vect{Z}}> M$:
$$
\pen[S,\alpha]{\Vect{X}}{\Vect{Z}}+\pen[S,\alpha]{\Vect{Y}}{\Vect{Z}}-\ent{\Vect{Z}}=(2\alpha-1)\ent{\Vect{Z}}+(1-\alpha)(m+M)\geq \alpha +(1-\alpha)M=\pen[S,\alpha]{\Vect{X}}{\Vect{Y}}.
$$
\item $m\leq \ent{\Vect{Z}} \leq M$:
$$
\pen[S,\alpha]{\Vect{X}}{\Vect{Z}}+\pen[S,\alpha]{\Vect{Y}}{\Vect{Z}}-\ent{\Vect{Z}}= \alpha m + (1-\alpha)M=\pen[S,\alpha]{\Vect{X}}{\Vect{Y}}.
$$
\end{itemize}

$(iv)$ trivial.
\end{proof}\newline

We assert that the measures $\IBempty[R,\alpha]$, \IBempty[P,\alpha]
and \IBempty[D,\alpha] (and so $\NIBempty[R,\alpha]$,
\NIBempty[P,\alpha] and \NIBempty[D,\alpha]) do not satisfy the
condition~(\ref{hyp1inegTri}). Consider for
example~\IBempty[D,\alpha]. Let us choose $\Vect{X},\Vect{Y}$ and
$\Vect{Z}$ such that $\ent{\Vect{Z}}>{\max\left(\ent{\Vect{X}} , \ent{\Vect{Y}} \right)}$ and such that
$\ent{\Vect{Z}}=\frac{1+\alpha}{\alpha}
\ent{\Vect{X}}=\frac{1+\alpha}{\alpha}\ent{\Vect{Y}}$. This leads to
\begin{eqnarray*}
\pen[D,\alpha]{\Vect{X}}{\Vect{Z}} + \pen[D,\alpha]{\Vect{Y}}{\Vect{Z}} - \ent{\Vect{Z}}  &=& \ent{\Vect{Z}} \left( \frac{\ent{\Vect{X}}}{\alpha
\ent{\Vect{X}}+ (1-\alpha)\ent{\Vect{Z}} } + \frac{\ent{\Vect{Y}}}{\alpha \ent{\Vect{Y}}+
(1-\alpha) \ent{\Vect{Z}} } -1 \right) \\ &=& 0 <\pen[D,\alpha]{\Vect{X}}{\Vect{Y}},
\end{eqnarray*}
which is in contradiction with (\ref{hyp1inegTri}).

Concerning these divergences and the measures \IBempty[S,\alpha]
(for $\alpha>\frac12$) and \NIBempty[S,\alpha], we do not know if
they satisfy a triangular inequality but our tool cannot be applied
to prove it. We propose to weaken the property \textbf{[P6]} in
the following way in order to obtain more results. An
\ibd or \nibd satisfies

\textbf{[P6bis($\Upsilon$,$c$)]} There exists $c\geq 1$ such that for all $(\Vect{X},\Vect{Y}),(\Vect{Y},\Vect{Z}),(\Vect{X},\Vect{Z}) \in
\Upsilon$
$$
\IB{\Vect{X}}{\Vect{Y}} \leq c \times \left( \IB{\Vect{X}}{\Vect{Z}}+ \IB{\Vect{Y}}{\Vect{Z}} \right).
$$
Property \textbf{[P6]} is then equivalent to \textbf{[P6bis($\Gamma^2$,$1$)]}
 and we already know that \DEempty, \dEempty,
\DIempty, \dIempty and \DIBempty[S,\alpha] (for $\alpha \leq
\frac12$) satisfy \textbf{[P6bis($\Gamma^2$,$1$)]}. When $\Upsilon \subsetneq \Gamma^2$ the property \textbf{[P6bis]} is in some sense local whereas it is global (as a classical triangular inequality) when $\Upsilon=\Gamma^2$.

Let us notice that if an \ibd (or \nibd) satisfies
\textbf{[P3bis($\Upsilon,k_1,k_2$)]}, then  \textbf{[P6bis($\Upsilon$,$\frac{k_2}{k_1}$)]} is satisfied since
$$
\IB{\Vect{X}}{\Vect{Y}} \leq k_2 \DIB[I]{\Vect{X}}{\Vect{Y}} \leq k_2 \left(
\DIB[I]{\Vect{X}}{\Vect{Z}}+\DIB[I]{\Vect{Y}}{\Vect{Z}}\right) \leq \frac{k_2}{k_1}
\left(\IB{\Vect{X}}{\Vect{Z}}+\IB{\Vect{Y}}{\Vect{Z}} \right).
$$
We then inherit a lot of results from Proposition~\ref{prop3-P3bis} related to our examples. In particular $\IBempty[\bullet,\alpha]$ and $\NIBempty[\bullet,\alpha]$ (where $\bullet$ stands for $S,R,P$ and $D$)  both satisfy \textbf{[P6bis($\Upsilon_\Theta$,$\frac{1}{k^a_1}$)]}, \textbf{[P6bis($\Gamma_\Theta^2$,$\frac{1}{k^b_1}$)]} and \textbf{[P6bis($\Gamma_\Theta^2$,$\frac{1}{k^c_1}$)]}.

In the rest of this section, we attempt to ensure the global property
\textbf{[P6bis($\Gamma^2$,$c$)]}.
From Proposition~\ref{prop3-P3bis} (with $\Theta=\RR^+$), we assert that the divergences \IBempty[S,\alpha] (when $\alpha>\frac12$) and \NIBempty[S,\alpha] (resp. \IBempty[R,\alpha] and \NIBempty[R,\alpha]) satisfy \textbf{[P6bis($\Gamma^2$,$\frac{1}{1-\alpha}$)]} (resp. \textbf{[P6bis($\Gamma^2$,$\frac{1}{(1-\alpha)^2}$)]}).

When $\alpha\leq \frac12$, we could improve the previous on \IBempty[R,\alpha] by proving that it satisfies \textbf{[P6bis($\Gamma^2$,$\frac1{\alpha^2+(1-\alpha)^2}$)]}, in the same spirit of the proof leading to~\textbf{[P3bis]}. Indeed,
\begin{eqnarray*}
\DIB[S,\alpha]{\Vect{X}}{\Vect{Y}}- \IB[R,\alpha]{\Vect{X}}{\Vect{Y}} &=& \alpha(1-\alpha)(m+M-2\sqrt{mM} \\
&\leq& 2\alpha(1-\alpha) \left( \DIB[S,\alpha]{\Vect{X}}{\Vect{Y}} + \IM{\Vect{X}}{\Vect{Y}} -\sqrt{mM}\right) \\
&\leq& 2\alpha(1-\alpha) \DIB[S,\alpha]{\Vect{X}}{\Vect{Y}}
\end{eqnarray*}
which leads to  $\IB[R,\alpha]{\Vect{X}}{\Vect{Y}} \geq (\alpha^2+(1-\alpha)^2) \DIB[S,\alpha]{\Vect{X}}{\Vect{Y}}$. Finally, let us notice that
$$\IB[R,\alpha]{\Vect{X}}{\Vect{Y}} \leq \DIB[S,\alpha]{\Vect{X}}{\Vect{Y}} \leq \DIB[S,\alpha]{\Vect{X}}{\Vect{Z}} + \DIB[S,\alpha]{\Vect{Y}}{\Vect{Z}} \leq
\frac1{\alpha^2+(1-\alpha)^2} \left(\IB[R,\alpha]{\Vect{X}}{\Vect{Z}}+ \IB[R,\alpha]{\Vect{Y}}{\Vect{Z}}\right).
$$

\noindent We now give a further and general result allowing us, in particular, to improve \textbf{[P6bis($\Gamma^2,\frac1{1-\alpha}$)]} for \IBempty[S,\alpha] when $\alpha>\frac12$.

\begin{proposition} \label{prop-inegGenInfoMes}
Let us consider the following assumptions on a \compl[:] there
exists a constant $c \geq 1$ such that
\begin{eqnarray}
c \;\pen{\Vect{X}}{\Vect{Z}} + c \;\pen{\Vect{Y}}{\Vect{Z}}  -
\ent{\Vect{Z}} - (c-1)\left(
\IM{\Vect{X}}{\Vect{Z}}+\IM{\Vect{Y}}{\Vect{Z}} \right) &\geq&
\pen{\Vect{X}}{\Vect{Y}} \label{hypInegGenIB} \\ c \;
\pen{\Vect{X}}{\Vect{Z}} + c \; \pen{\Vect{Y}}{\Vect{Z}}  -
\ent{\Vect{Z}} - (c-1)\left(
\IM{\Vect{X}}{\Vect{Z}}+\IM{\Vect{Y}}{\Vect{Z}} \right) &\geq& \max
\left(
\pen{\Vect{X}}{\Vect{Y}},\pen{\Vect{X}}{\Vect{Z}},\pen{\Vect{Y}}{\Vect{Z}}\right).
\label{hypInegGen}
\end{eqnarray}
If an \ibd satisfies (\ref{hypInegGenIB}) or a \nibd
satisfies~(\ref{hypInegGen}), then they
satisfy~\textbf{[P6bis($\Gamma^2$,$c$)]}.
\end{proposition}

\begin{proof}
Let us introduce
$$
A =-\left(\pen{\Vect{X}}{\Vect{Y}} -\IM{\Vect{X}}{\Vect{Y}}\right) +
c\times( \pen{\Vect{X}}{\Vect{Z}} - \IM{\Vect{X}}{\Vect{Z}} ) +
c\times( \pen{\Vect{Y}}{\Vect{Z}} - \IM{\Vect{Y}}{\Vect{Z}}).
$$From~(\ref{eq-resIMTri}) and (\ref{hypInegGenIB}), one may assert that
$$
A \geq c \;\pen{\Vect{X}}{\Vect{Z}} + c \; \pen{\Vect{Y}}{\Vect{Z}}
- \pen{\Vect{X}}{\Vect{Y}} - \ent{\Vect{Z}} - (c-1)\left(
\IM{\Vect{X}}{\Vect{Z}}+\IM{\Vect{Y}}{\Vect{Z}} \right)\geq 0,
$$
which implies that the result is valid for $\IBempty$. Now, from~(\ref{hypInegGen}) one can write
$$A+\pen{\Vect{X}}{\Vect{Y}} \geq \max \left( \pen{\Vect{X}}{\Vect{Z}},\pen{\Vect{Y}}{\Vect{Z}}\right)
$$
which leads to
$$
\NIB{\Vect{X}}{\Vect{Y}} \leq \frac{c\times
\left(\pen{\Vect{X}}{\Vect{Z}} - \IM{\Vect{X}}{\Vect{Z}} \right) +
c\times \left(\pen{\Vect{Y}}{\Vect{Z}} - \IM{\Vect{Y}}{\Vect{Z}}
\right)}{ \max \left(
\pen{\Vect{X}}{\Vect{Z}},\pen{\Vect{Y}}{\Vect{Z}}\right)} \leq c
\times \NIB{\Vect{X}}{\Vect{Z}} + c \times \NIB{\Vect{Y}}{\Vect{Z}}.
$$
\end{proof}

\begin{corollary} \label{cor-ineghE} ${ }$

The measures \IBempty[S,\alpha] for $\alpha>\frac12$ satisfy \textbf{[P6bis($\Gamma^2$,$\frac{\alpha}{1-\alpha}$)]}

\end{corollary}

\begin{proof}
Let us concentrate on \IBempty[S,\alpha] for $\alpha>\frac12$.
Let $A=c \pen[S,\alpha]{\Vect{X}}{\Vect{Z}}+c \pen[S,\alpha]{\Vect{Y}}{\Vect{Z}} -\ent{\Vect{Z}}
-(c-1)\left( \IM{\Vect{X}}{\Vect{Z}}+\IM{\Vect{Y}}{\Vect{Z}} \right)$. Without loss of generality,
we assume $\ent{x} \leq \ent{\Vect{Y}}$. We distinguish three cases:
\begin{itemize}
\item $\ent{\Vect{Z}} \leq \ent{\Vect{X}}\leq \ent{\Vect{Y}}$: we have
$$
A \geq  c(1-\alpha)\ent{\Vect{X}} + (1-\alpha)\ent{\Vect{Y}}+
(c\alpha+\alpha-1)\ent{\Vect{Z}}-(c-1)\IM{\Vect{X}}{\Vect{Z}}.
$$
Then,
$$
A-\pen[S,\alpha]{\Vect{X}}{\Vect{Y}}\geq (c(1-\alpha)-\alpha)\ent{\Vect{X}}+
(c\alpha+\alpha-1)\ent{\Vect{Z}} -(c-1)\IM{\Vect{X}}{\Vect{Z}} \geq (c-1)\left(
\ent{\Vect{Z}}-\IM{\Vect{X}}{\Vect{Z}} \right) \geq 0,
$$
as soon as $c\geq \frac{\alpha}{1-\alpha}$.
\item $\ent{\Vect{X}} \leq \ent{\Vect{Y}} \leq \ent{\Vect{Z}}$: we have
$$
A\geq \alpha\ent{\Vect{X}} + c\alpha \ent{\Vect{Y}}+
\left((1-\alpha)+c(1-\alpha)-1 \right)\ent{\Vect{Z}}-(c-1)\IM{\Vect{Y}}{\Vect{Z}}.
$$
Then,
$$
A-\pen[S,\alpha]{\Vect{X}}{\Vect{Y}}\geq (c\alpha-(1-\alpha))\ent{\Vect{Y}} +
((1-\alpha)+c(1-\alpha)-1)\ent{\Vect{Z}} -(c-1)\IM{\Vect{Y}}{\Vect{Z}} \geq (c-1)\left(
\ent{\Vect{Y}}-\IM{\Vect{Y}}{\Vect{Z}}\right)\geq 0,
$$
as soon as $c\geq \frac{\alpha}{1-\alpha}$.
\item $\ent{\Vect{X}} <\ent{\Vect{Z}}< \ent{\Vect{Y}}$: we have
$$
A \geq c\alpha \ent{\Vect{X}} + (1-\alpha)\ent{\Vect{Y}} + (c(1-\alpha)\ent{\Vect{Z}} +
\alpha-1) - (c-1)\IM{\Vect{X}}{\Vect{Z}}.
$$
Then,
$$
A-\pen[S,\alpha]{\Vect{X}}{\Vect{Y}}\geq (c-1)\alpha\ent{\Vect{X}}+
(c-1)(1-\alpha)\ent{\Vect{Z}} -\IM{\Vect{X}}{\Vect{Z}} \geq 0.
$$
\end{itemize}
Hence, we obtain for $c=\frac{\alpha}{1-\alpha}$, $A-\pen[S,\alpha]{\Vect{X}}{\Vect{Y}}\geq 0$. \\
\end{proof}

\begin{remark} \label{rem-hE}
The tool presented in Proposition~\ref{prop-inegGenInfoMes} cannot
be applied to the \ibd $\IBempty[D,\alpha]$ and the \nibd
$\NIBempty[D,\alpha]$. Indeed, let us give some $c\geq1$ and let us
consider the quantity
$$
A=c \; \pen[D,\alpha]{\Vect{X}}{\Vect{Z}} +
c\;\pen[D,\alpha]{\Vect{Y}}{\Vect{Z}}-\ent{\Vect{Z}} -(c-1) \left(
\IM{\Vect{X}}{\Vect{Z}} + \IM{\Vect{Y}}{\Vect{Z}} \right).
$$
In fact, one can always find $\Vect{X},\Vect{Y},\Vect{Z}$ such that
for all $c\geq 1$, the quantity $A$ is negative. Indeed, let us
choose $\Vect{Z}$ independent of $\Vect{X}$ and $\Vect{Y}$ and such
that
$\alpha\ent{\Vect{Z}}+(1-\alpha)\ent{\Vect{X}}=3c\ent{\Vect{X}}$ and
$\alpha\ent{\Vect{Z}}+(1-\alpha)\ent{\Vect{Y}}=3c\ent{\Vect{Y}}$.
Then, it is easy to see that
$A=\ent{\Vect{Z}}\left(\frac13+\frac13-1\right) < 0$. In the same
manner, the tool is inapplicable to the \ibd $\IBempty[P,\alpha]$
and the \nibd $\NIBempty[P,\alpha]$. Indeed, let us give
$\Vect{Z}$ independent of $\Vect{X}$ and $\Vect{Y}$ and such that
$\ent{\Vect{X}}=\ent{\Vect{Y}} =\left(\frac{1}{3c}\right)^{1/\alpha} \ent{\Vect{Z}}$, then
$$
A= c \; \pen[P,\alpha]{\Vect{X}}{\Vect{Z}} + c\;\pen[P,\alpha]{\Vect{Y}}{\Vect{Z}}-\ent{\Vect{Z}}
-(c-1) \left( \IM{\Vect{X}}{\Vect{Z}} + \IM{\Vect{Y}}{\Vect{Z}} \right) = -\frac13 \ent{\Vect{Z}}<0.
$$
\end{remark}

The following result is an extension of
Proposition~\ref{prop-inegGenInfoMes} well-suited to be applied to
$\NIBempty[D,\alpha]$.

\begin{proposition} \label{prop-tool4hE}
Let us assume that there exists two positive integer $I$ and $J$ such that a \nibd $\NIB{\Vect{X}}{\Vect{Y}}$
can be expressed as:
\[
\NIB{\Vect{X}}{\Vect{Y}}=\sum_{i=1}^I \frac{S^{(i)}_{\Vect{X},\Vect{Y}}}{U^{(i)}_{\Vect{X},\Vect{Y}}}=\sum_{j=1}^J \alpha^{(j)}\left(1-\frac{\IM{\Vect{X}}{\Vect{Y}}}{\pen[(j)]{\Vect{X}}{\Vect{Y}}}\right)
\]
where $\left(\alpha^{(j)}\right)_{j=1,\cdots,J}$ is some vector of
probability weights. By denoting $S_{\Vect{X},\Vect{Y}}=\sum_{i=1}^I
S^{(i)}_{\Vect{X},\Vect{Y}}$ and
$U_{\Vect{X},\Vect{Y}}=\max_{i=1,\cdots,I}
U^{(i)}_{\Vect{X},\Vect{Y}}$, if there exists some real number
$c\geq 1$ such that for any $j=1,\cdots,J$ the following assumptions
are satisfied:
\begin{itemize}
\item[(i)] $A^{(j)}=\IM{\Vect{X}}{\Vect{Y}}-\pen[(j)]{\Vect{X}}{\Vect{Y}} +c\left(S_{\Vect{X},\Vect{Z}}+S_{\Vect{Z},\Vect{Y}}\right)\geq 0$.
\item[(ii)] $A^{(j)}+\pen[(j)]{\Vect{X}}{\Vect{Y}} \geq \max(U_{\Vect{X},\Vect{Z}},U_{\Vect{Z},\Vect{Y}})$.
\end{itemize}
then \NIBempty satisfies \textbf{[P6bis($\Gamma^2$,$c$)]}.
\end{proposition}

\begin{proof}
Using assumptions $(i)$ and $(ii)$, one can prove that for all $j=1,\ldots,J$
\begin{eqnarray}
1-\frac{\IM{\Vect{X}}{\Vect{Y}}}{\penk{j}{\Vect{X}}{\Vect{Y}}} &\leq& 1-
\frac{ \IM{\Vect{X}}{\Vect{Y}}}{\penk{j}{\Vect{X}}{\Vect{Y}}+A^{(j)}} \leq
c\frac{S_{\Vect{X},\Vect{Z}}+S_{\Vect{Y},\Vect{Z}}}{\max(U_{\Vect{X},\Vect{Z}},U_{\Vect{Z},\Vect{Y}})} \nonumber \\
&\leq & c \times( \NIB{\Vect{X}}{\Vect{Z}}+ \NIB{\Vect{Y}}{\Vect{Z}}). \nonumber
\end{eqnarray}
It follows that
$$
\NIB{\Vect{X}}{\Vect{Y}}=\sum_{j=1}^J \alpha^{(j)}\left(1-\frac{\IM{\Vect{X}}{\Vect{Y}}}{\pen[(j)]{\Vect{X}}{\Vect{Y}}}\right) \leq \sum_{j=1}^J \alpha^{(j)} \times c \times( \NIB{\Vect{X}}{\Vect{Z}}+ \NIB{\Vect{Y}}{\Vect{Z}}) = c \times( \NIB{\Vect{X}}{\Vect{Z}}+ \NIB{\Vect{Y}}{\Vect{Z}}).
$$
\end{proof}

\begin{corollary}
The measure $\NIBempty[D,\alpha]$ satisfies~\textbf{[P6bis($\Gamma^2$,$\frac1{{\alpha^\wedge}}$)]}.
\end{corollary}

\begin{proof}
We have
\begin{eqnarray*}
\NIB[D,\alpha]{\Vect{X}}{\Vect{Y}}&=& \alpha \min \left(
\frac{\ent{\Given{\Vect{X}}{\Vect{Y}}}}{\ent{\Vect{X}}} ,
\frac{\ent{\Given{\Vect{Y}}{\Vect{X}}}}{\ent{\Vect{Y}}}\right) +
(1-\alpha) \max \left(
\frac{\ent{\Given{\Vect{X}}{\Vect{Y}}}}{\ent{\Vect{X}}} ,
\frac{\ent{\Given{\Vect{Y}}{\Vect{X}}}}{\ent{\Vect{Y}}}\right) \\
&=& \alpha \frac{ \min \left(\ent{\Given{\Vect{X}}{\Vect{Y}}}, \ent{\Given{\Vect{Y}}{\Vect{X}}}\right)}{{\min\left(\ent{\Vect{X}} , \ent{\Vect{Y}} \right)}} + (1-\alpha)\frac{ \max \left(\ent{\Given{\Vect{X}}{\Vect{Y}}}, \ent{\Given{\Vect{Y}}{\Vect{X}}}\right)}{{\max\left(\ent{\Vect{X}} , \ent{\Vect{Y}} \right)}} \\
&=& \alpha\left( 1-\frac{\IM{\Vect{X}}{\Vect{Y}}}{{\min\left(\ent{\Vect{X}} , \ent{\Vect{Y}} \right)}} \right)
+(1-\alpha)\left( 1-\frac{\IM{\Vect{X}}{\Vect{Y}}}{{\max\left(\ent{\Vect{X}} , \ent{\Vect{Y}} \right)}} \right).
\end{eqnarray*}
By identification with notation introduced in
Proposition~\ref{prop-tool4hE}, we have $I=J=2$,
$S^{(1)}_{\Vect{X},\Vect{Y}}=\alpha \min \left(\ent{\Given{\Vect{X}}{\Vect{Y}}}, \ent{\Given{\Vect{Y}}{\Vect{X}}}\right) $,
$S^{(2)}_{\Vect{X},\Vect{Y}}=(1-\alpha)  \max\left(\ent{\Given{\Vect{X}}{\Vect{Y}}}, \ent{\Given{\Vect{Y}}{\Vect{X}}}\right) $,
$U^{(1)}_{\Vect{X},\Vect{Y}}={\min\left(\ent{\Vect{X}} , \ent{\Vect{Y}} \right)}$,
$U^{(2)}_{\Vect{X},\Vect{Y}}={\max\left(\ent{\Vect{X}} , \ent{\Vect{Y}} \right)}$,
$\penk{1}{\Vect{X}}{\Vect{Y}}={\min\left(\ent{\Vect{X}} , \ent{\Vect{Y}} \right)}$ and
$\penk{2}{\Vect{X}}{\Vect{Y}}={\max\left(\ent{\Vect{X}} , \ent{\Vect{Y}} \right)}$. Let us fix $c$ to the value
$\frac1{{\alpha^\wedge}}$. We have
\begin{eqnarray*}
A^{(1)} &=& \IM{\Vect{X}}{\Vect{Y}} -{\min\left(\ent{\Vect{X}} , \ent{\Vect{Y}} \right)} + \frac{1}{{\alpha^\wedge}}
\left( \alpha \min \left( \ent{\Given{\Vect{X}}{\Vect{Z}}} ,
\ent{\Given{\Vect{Z}}{\Vect{X}}} \right) + (1-\alpha) \max \left(
\ent{\Given{\Vect{X}}{\Vect{Z}}} , \ent{\Given{\Vect{Z}}{\Vect{X}}}
\right) \right. \\ && \left. + \alpha \min \left(
\ent{\Given{\Vect{Y}}{\Vect{Z}}} , \ent{\Given{\Vect{Z}}{\Vect{Y}}}
\right) + (1-\alpha) \max \left( \ent{\Given{\Vect{Y}}{\Vect{Z}}} ,
\ent{\Given{\Vect{Z}}{\Vect{Y}}} \right) \right)
\end{eqnarray*}
Clearly from~(\ref{eq-resEntJoint})
\begin{eqnarray}
A^{(1)}  &\geq& {\max\left(\ent{\Vect{X}} , \ent{\Vect{Y}} \right)} - \ent{\Vect{X},\Vect{Y}} + 2 \ent{\Vect{X},\Vect{Z}} + 2 \ent{\Vect{Y},\Vect{Z}} -\ent{\Vect{X}} -\ent{\Vect{Y}} -2\ent{\Vect{Z}} \nonumber\\
&\geq & \ent{\Vect{X},\Vect{Z}} + \ent{\Vect{Y},\Vect{Z}}-{\min\left(\ent{\Vect{X}} , \ent{\Vect{Y}} \right)} -\ent{\Vect{Z}}\geq 0.\nonumber
\end{eqnarray}
And one also has
$$
{\min\left(\ent{\Vect{X}} , \ent{\Vect{Y}} \right)}+A^{(1)} \geq \ent{\Vect{X},\Vect{Z}} + \ent{\Vect{Y},\Vect{Z}} -\ent{\Vect{Z}} \geq \max\left(\ent{\Vect{X}},\ent{\Vect{Y}},\ent{\Vect{Z}} \right)= \max\left(U_{\Vect{X},\Vect{Z}}, U_{\Vect{Y},\Vect{Z}}\right).
$$
It follows that $A^{(1)}$ fullfills conditions $(i)$ and $(ii)$ of
Proposition~\ref{prop-tool4hE} with $c=\frac1{{\alpha^\wedge}}$. The proof
is strictly similar for $A^{(2)}$.
\end{proof}

\section{Prediction framework} \label{sec-pred}

We pay attention on properties related to the prediction of some fixed random vector~$\Vect{Y}$.

\subsection{Prediction framework}

Recall that our purpose is to find the random vector $\Vect{X}$ that minimizes  $\IB{\Vect{Y}}{\Vect{X}}$ (resp. $\NIB{\Vect{Y}}{\Vect{X}}$) which combines a complexity term $\pen{\Vect{X}}{\Vect{Y}}$ (to minimize) and an information term $\IM{\Vect{X}}{\Vect{Y}}$ (to maximize). Let us imagine that we already get some $\Vect{X}_1$ and its associated measure $\IB{\Vect{Y}}{\Vect{X}_1}$ (resp. $\NIB{\Vect{Y}}{\Vect{X}_1}$). After evaluating $\IB{\Vect{Y}}{\Vect{X}_2}$ (resp. $\NIB{\Vect{Y}}{\Vect{X}_2}$), we may be interested in describing the conditions under which $\Vect{X}_2$ is better or worse than $\Vect{X}_1$:

\begin{proposition} \label{prop-compNIB}
Two situations may occur

\noindent Case 1: we choose $\Vect{X}_2$ instead of $\Vect{X}_1$  when
\begin{eqnarray} 
\IB{\Vect{Y}}{\Vect{X}_2} < \IB{\Vect{Y}}{\Vect{X}_1} &\Longleftrightarrow&
  \pen{\Vect{Y}}{\Vect{X}_2}-\pen{\Vect{Y}}{\Vect{X}_1} < \IM{\Vect{Y}}{\Vect{X}_2}- \IM{\Vect{Y}}{\Vect{X}_1}\label{eq-interetIB} \\
\NIB{\Vect{Y}}{\Vect{X}_2} < \NIB{\Vect{Y}}{\Vect{X}_1} &\Longleftrightarrow&   \frac{
\pen{\Vect{Y}}{\Vect{X}_2}- \pen{\Vect{Y}}{\Vect{X}_1}}{\pen{\Vect{Y}}{\Vect{X}_1}} <        \frac{\IM{Y}{\Vect{X}_2}- \IM{Y}{\Vect{X}_1}}{\IM{Y}{\Vect{X}_1}} \label{eq-interetNIB}
\end{eqnarray}
Case 2: we keep $\Vect{X}_1$ and reject $\Vect{X}_2$ when
\begin{eqnarray} 
\IB{\Vect{Y}}{\Vect{X}_2} \geq \IB{\Vect{Y}}{\Vect{X}_1} &\Longleftrightarrow&
  \pen{\Vect{Y}}{\Vect{X}_2}-\pen{\Vect{Y}}{\Vect{X}_1} \geq \IM{\Vect{Y}}{\Vect{X}_2}- \IM{\Vect{Y}}{\Vect{X}_1}\label{eq-interetIB} \\
\NIB{\Vect{Y}}{\Vect{X}_2} \geq \NIB{\Vect{Y}}{\Vect{X}_1} &\Longleftrightarrow&   \frac{
\pen{\Vect{Y}}{\Vect{X}_2}- \pen{\Vect{Y}}{\Vect{X}_1}}{\pen{\Vect{Y}}{\Vect{X}_1}} \geq        \frac{\IM{Y}{\Vect{X}_2}- \IM{Y}{\Vect{X}_1}}{\IM{Y}{\Vect{X}_1}} \label{eq-interetNIB}
\end{eqnarray}
\end{proposition}

This result implies automatically that the properties \textbf{[P8]} and \textbf{[P9]}
are satisfied. Let us comment more precisely the previous proposition:
\begin{itemize}
\item Case 1 holds when
\begin{enumerate}
\item $\Vect{X}_2$ is simpler than $\Vect{X}_1$ (i.e. $\pen{\Vect{Y}}{\Vect{X}_2}-\pen{\Vect{Y}}{\Vect{X}_1}< 0$) and $\Vect{X}_2$ is at least as informative as $\Vect{X}_1$ (i.e.  $\IM{\Vect{Y}}{\Vect{X}_2}-\IM{\Vect{Y}}{\Vect{X}_1}\geq 0$).
\item $\Vect{X}_2$ and $\Vect{X}_1$ have the same complexity (i.e. $\pen{\Vect{Y}}{\Vect{X}_2}-\pen{\Vect{Y}}{\Vect{X}_1}= 0$) and $\Vect{X}_2$ is more informative than $\Vect{X}_1$ (i.e.  $\IM{\Vect{Y}}{\Vect{X}_2}-\IM{\Vect{Y}}{\Vect{X}_1}>0$).
\item $\Vect{X}_2$ is simpler and less informative than $\Vect{X}_1$ and such that the absolute (resp. relative) excess of complexity is lower than the absolute (resp. relative) gain of information that is $\pen{\Vect{Y}}{\Vect{X}_2}-\pen{\Vect{Y}}{\Vect{X}_1}<\IM{\Vect{Y}}{\Vect{X}_2}-\IM{\Vect{Y}}{\Vect{X}_1}<0 $ (resp. $\frac{\pen{\Vect{Y}}{\Vect{X}_2}-\pen{\Vect{Y}}{\Vect{X}_1}}{\pen{\Vect{Y}}{\Vect{X}_1} }< \frac{\IM{\Vect{Y}}{\Vect{X}_2}-\IM{\Vect{Y}}{\Vect{X}_1}}{ \IM{\Vect{Y}}{\Vect{X}_1}}<0$).
\item $\Vect{X}_2$ is more complex and more informative than $\Vect{X}_1$ and such that the absolute (resp. relative) excess of complexity is lower than the absolute (resp. relative) gain of information that is $0<\pen{\Vect{Y}}{\Vect{X}_2}-\pen{\Vect{Y}}{\Vect{X}_1}<\IM{\Vect{Y}}{\Vect{X}_2}-\IM{\Vect{Y}}{\Vect{X}_1}$ (resp. $0<\frac{\pen{\Vect{Y}}{\Vect{X}_2}-\pen{\Vect{Y}}{\Vect{X}_1}}{\pen{\Vect{Y}}{\Vect{X}_1} }< \frac{\IM{\Vect{Y}}{\Vect{X}_2}-\IM{\Vect{Y}}{\Vect{X}_1}}{ \IM{\Vect{Y}}{\Vect{X}_1}}$).
\end{enumerate}
\item Case 2 holds when
\begin{enumerate}
\item $\Vect{X}_2$ is at least as complex as $\Vect{X}_1$ (i.e. $\pen{\Vect{Y}}{\Vect{X}_2}-\pen{\Vect{Y}}{\Vect{X}_1}\geq 0$) and $\Vect{X}_2$ is at most as informative as $\Vect{X}_1$ (i.e.  $\IM{\Vect{Y}}{\Vect{X}_2}-\IM{\Vect{Y}}{\Vect{X}_1}\leq 0$).
\item $\Vect{X}_2$ is simpler and less informative than $\Vect{X}_1$, and such that the 
absolute (resp. relative) excess of complexity is greater than or equal to the absolute (resp. relative) gain of information that is $0>\pen{\Vect{Y}}{\Vect{X}_2}-\pen{\Vect{Y}}{\Vect{X}_1}\geq \IM{\Vect{Y}}{\Vect{X}_2}-\IM{\Vect{Y}}{\Vect{X}_1}$ (resp. $0>\frac{\pen{\Vect{Y}}{\Vect{X}_2}-\pen{\Vect{Y}}{\Vect{X}_1}}{\pen{\Vect{Y}}{\Vect{X}_1} }\geq \frac{\IM{\Vect{Y}}{\Vect{X}_2}-\IM{\Vect{Y}}{\Vect{X}_1}}{ \IM{\Vect{Y}}{\Vect{X}_1}}$).
\item $\Vect{X}_2$ is more complex and more informative than $\Vect{X}_1$, and such that the 
absolute (resp. relative) excess of complexity is greater than or equal to the absolute (resp. relative) gain of information that is $\pen{\Vect{Y}}{\Vect{X}_2}-\pen{\Vect{Y}}{\Vect{X}_1}\geq \IM{\Vect{Y}}{\Vect{X}_2}-\IM{\Vect{Y}}{\Vect{X}_1}>0$ (resp. $\frac{\pen{\Vect{Y}}{\Vect{X}_2}-\pen{\Vect{Y}}{\Vect{X}_1}}{\pen{\Vect{Y}}{\Vect{X}_1} }\geq \frac{\IM{\Vect{Y}}{\Vect{X}_2}-\IM{\Vect{Y}}{\Vect{X}_1}}{ \IM{\Vect{Y}}{\Vect{X}_1}}>0$).
\end{enumerate}
\end{itemize}

\begin{proposition} \label{prop-Pcompl}
Any \compl \penempty[\alpha] of the form~(\ref{eq-NIBg}) satisfies \textbf{[P10]}.
\end{proposition}

\begin{proof}
Without loss of generality the function $g(\cdot)$ defining \penempty[\alpha] is assumed to be an increasing function. 
Hence, $\ent{\Vect{X}_2}\geq \ent{\Vect{X}_1}$ implies that $\pen[\alpha]{\Vect{Y}}{\Vect{X}_2}\geq  \pen[\alpha]{\Vect{Y}}{\Vect{X}_1}$. Now, let us assume $\pen[\alpha]{\Vect{Y}}{\Vect{X}_2}\geq  \pen[\alpha]{\Vect{Y}}{\Vect{X}_1}$. We assert by denoting $m_i=\min(\ent{\Vect{Y}},\ent{{\Vect{X}}_i})$ and $M_i=\max(\ent{\Vect{Y}},\ent{{\Vect{X}}_i})$ for $i=1,2$
\begin{eqnarray*}
\pen[\alpha]{\Vect{Y}}{\Vect{X}_2}\geq  \pen[\alpha]{\Vect{Y}}{\Vect{X}_1} &\Longleftrightarrow& g^{-1} \left( \alpha g(m_1)+ (1-\alpha)g(M_1)\right) \leq g^{-1} \left( \alpha g(m_2)+ (1-\alpha)g(M_2)\right) \\
&\Longleftrightarrow& \alpha \left(g(m_1)-g(m_2) \right) + (1-\alpha) \left(g(M_1)-g(M_2) \right) \leq 0
\end{eqnarray*}
Now, assume moreover that $\ent{\Vect{X}_1}> \ent{\Vect{X}_2}$, then the right-hand side is 
$$
\left\{ 
\begin{array}{lc}
=(1-\alpha)(g(\ent{\Vect{X}_1}-g(\ent{\Vect{X}_2}))>0 & \mbox{ if } \ent{\Vect{Y}}\leq \ent{\Vect{X}_2}< \ent{\Vect{X}_1} \\
> g(m_1)-g(m_2) =0 & \mbox{ if } \ent{\Vect{X}_2}< \ent{\Vect{Y}}< \ent{\Vect{X}_1} \\
= \alpha(g(\ent{\Vect{X}_2})-g(\ent{\Vect{X}_1}))>0 & \mbox{ if } \ent{\Vect{X}_2}< \ent{\Vect{X}_1}\leq \ent{\Vect{Y}} 
\end{array}
\right. .
$$
This leads to a contradiction which implies that $\ent{\Vect{X}_2}\geq \ent{\Vect{X}_1}$.
\end{proof}

\begin{remark} \label{rem-Pcompl}
The \compl[s] \penempty[E] and \penempty[I] do not satisfy the property \textbf{[P10]} in the general case. Indeed, there is no implication for \penempty[E] and one can only prove that $\ent{\Vect{X}_1} \geq \ent{\Vect{X}_2} \Rightarrow \pen[I]{\Vect{Y}}{\Vect{X}_1}\geq \pen[I]{\Vect{Y}}{\Vect{X}_2}$. However, one can point out that when $\IM{\Vect{Y}}{\Vect{X}_1}=\IM{\Vect{Y}}{\Vect{X}_2}$ then both \penempty[E] and \penempty[I] satisfy \textbf{[P10]}.
\end{remark}

More specifically, two frameworks may be of special interest:
\begin{itemize}
\item $\Vect{X}_2$ is as informative as $\Vect{X}_1$ (i.e. $\IM{\Vect{Y}}{\Vect{X}_1}=\IM{\Vect{Y}}{\Vect{X}_2}$):
we expect to select the random variable with the smallest entropy. This is effectively what happens when \textbf{[P10]} which is satisfied from Proposition~\ref{prop-Pcompl} and Remark~\ref{rem-Pcompl} (in this framework)

 $\penempty[\bullet]$ with $\bullet=I, S, R, P, D$ in the general case and for $\penempty[E]$ in this framework since $\ent{\Vect{Y},\Vect{X}_2}-\ent{\Vect{Y},\Vect{X}_1}= \ent{\Vect{X}_2}-\ent{\Vect{X}_1}$. 

\item ${\Vect{X}_1}=g({\Vect{X}_2})$ with g some surjective (but not injective) mapping: 
${\Vect{X}_2}$ is more complex than ${\Vect{X}_1}$ and ${\Vect{X}_2}$ is at least as informative as ${\Vect{X}_1}$.
Consequently, this case is not trivial since
both absolute (resp. relative) excess of complexity and absolute (resp. relative) gain of information are competing. Let us give two important examples of such a context. 
\begin{enumerate}
\item quantization problem: given a quantized version $\Vect{X}_1$ of some (continuous) random variable with its associated partition $\mathcal{A}_1$, the problem is to know whether some new quantized version $\Vect{X}_2$ with an associated partition $\mathcal{A}_2$ finer than $\mathcal{A}_1$ should be preferred to predict $\Vect{Y}$.
\item variables selection problem: suppose one wants to construct an ascending selection method. The vector $\Vect{X}_1$ could represent some selected set of covariables and $\Vect{X}_2=(\Vect{X}_1,\Vect{X}_2^\prime)$ a larger set of covariables. The aim is so to know if $\Vect{X}_2^\prime$ should be integrated to the selected set or not.
\end{enumerate}
Some simple algorithms of quantization and selection methods are proposed in~\cite{Robineau04} using these results.
\end{itemize}

\bigskip

\subsection{Around the redundancy of two random vectors $\Vect{X}_1$ and $\Vect{X}_2$}

In the future use of an \ibd or \nibd[,] one would expect that if two discrete-valued random vectors $\Vect{X}_1$ and $\Vect{X}_2$ have the same  (or almost the same) information with respect  to an \ibd or  \nibd[,] then both have the same effect on the prediction of another vector $\Vect{Y}$. This requirement, expressed by the property \textbf{[P11]}, could be used for example in a variables selection problem in the context of discrimination to detect redundant variables.

In order to make the property \textbf{[P11]} applicable for practical purpose, we may find interesting to have a bound of the difference $|\IB{\Vect{Y}}{\Vect{X}_1}-\IB{\Vect{Y}}{\Vect{X}_2}|$ (resp. $|\NIB{\Vect{Y}}{\Vect{X}_1}-\NIB{\Vect{Y}}{\Vect{X}_2}|$) expressed in terms of $\DIB[I]{\Vect{X}_1}{\Vect{X}_2}$ (resp. $\dIB[I]{\Vect{X}_1}{\Vect{X}_2}$). More precisely, the question may arise whether there exists a function $h(\cdot)$ satisfying  $h(x)\to 0$ as $x\to 0$ and such that $|\IB{\Vect{Y}}{\Vect{X}_1}-\IB{\Vect{Y}}{\Vect{X}_2}|\leq h(\DIB[I]{\Vect{X}_1}{\Vect{X}_2})$ (resp. $|\NIB{\Vect{Y}}{\Vect{X}_1}-\NIB{\Vect{Y}}{\Vect{X}_2}|\leq h(\dIB[I]{\Vect{X}_1}{\Vect{X}_2})$). Here, according to our examples, we only concentrate ourself on linear function $h(\cdot)$.

We then propose to translate the property \textbf{[P11]} on an \ibd \IBempty (resp. a \nibd \NIBempty) by:

\textbf{[P11bis($\Upsilon,k$)]} there exists some positive constant $k$ such that for all $\sloppy{({\Vect{X}_1},{\Vect{X}_2}) \in \Upsilon\subset\Gamma^2}$ such that
\begin{equation} \label{eq-P10bis}
| \IB{\Vect{Y}}{\Vect{X}_1} - \IB{\Vect{Y}}{\Vect{X}_2} | \leq k \;  \DIB[I]{\Vect{X}_1}{\Vect{X}_2} 
\end{equation}

As a first answer, let us precise that if the \ibd (resp. \nibd[]) satisfies a triangular inequality \textbf{[P6bis($\Gamma^2,1$)]} and \textbf{[P3bis($\Upsilon,k_1,k_2$)]} then it satisfies \textbf{[P11bis($\Upsilon,k_2$)]} due to the equivalent expression of the triangular inequality as
$$
| \DIB{\Vect{Y}}{\Vect{X}_1} - \DIB{\Vect{Y}}{\Vect{X}_2} | \leq \DIB{\Vect{X}_1}{\Vect{X}_2}
\quad \mbox{ (resp. } \quad
| \dIB{\Vect{Y}}{\Vect{X}_1} - \dIB{\Vect{Y}}{\Vect{X}_2} | \leq \dIB{\Vect{X}_1}{\Vect{X}_2}\mbox{)}.
$$

A priori, if an \ibd or \nibd only satisfies \textbf{[P6bis($\Gamma^2,c$)]} with some $c> 1$, then this property does no more seem to be true: indeed, for all $\Vect{Y},\Vect{X}_1$ and $\Vect{X}_2$, one may prove for an \ibd by instance that
$$
| \IB{\Vect{Y}}{\Vect{X}_1} - \IB{\Vect{Y}}{\Vect{X}_2} | \leq c\times\IB{\Vect{X}_1}{\Vect{X}_2} + (c-1)\min\left(\IB{\Vect{Y}}{\Vect{X}_1},\IB{\Vect{Y}}{\Vect{X}_2}\right)\nleq c\times\IB{\Vect{X}_1}{\Vect{X}_2}.
$$
 Actually, this apparent disappointing result only expresses that a ``redundancy'' property cannot (always) be derived from a triangular's type inequality. \\

The following proposition gives some sufficient conditions required on some \compl ensuring that the associated \IBempty and \NIBempty satisfies the property \textbf{[P11bis]}

\begin{proposition} \label{prop-redundant}
$(i)$ Assume there exists some positive constant $\kappa_1$ such that the \compl of an \ibd satisfies for all $({\Vect{X}_1},{\Vect{X}_2}) \in \Upsilon$
\begin{equation} \label{hyp1Pen}
\left| \pen{\Vect{Y}}{\Vect{X}_1}  - \pen{\Vect{Y}}{\Vect{X}_2}  \right| \leq \kappa_1 \left| \ent{\Vect{X}_1}-\ent{\Vect{X}_2}\right|,
\end{equation}
then \IBempty satisfies \textbf{[P11bis$(\Upsilon,1+\kappa_1)$]}

\noindent $(ii)$ If in addition, there exists some positive constant $\kappa_2$ such that for all $({\Vect{X}_1},{\Vect{X}_2}) \in \Upsilon$
\begin{equation} \label{hyp2Pen}
\max\left( \pen{\Vect{Y}}{\Vect{X}_1}, \pen{\Vect{Y}}{\Vect{X}_2}\right)   \geq \kappa_2 \times \pen[I]{\Vect{X}_1}{\Vect{X}_2} 
\end{equation}
then the associated \nibd satisfies \textbf{[P11bis$\left(\Upsilon,\frac{1+\kappa_1}{\kappa_2}\right)$]}
\end{proposition}

\begin{proof}
$(i)$ Let us start to write
\begin{equation} \label{eq1-proofPropDiff}
\left| \IB{\Vect{Y}}{\Vect{X}_1}- \IB{\Vect{Y}}{\Vect{X}_2} \right| \leq \left| \IM{\Vect{Y}}{\Vect{X}_1}- \IM{\Vect{Y}}{\Vect{X}_2} \right|+ \left| \pen{\Vect{Y}}{\Vect{X}_1}- \pen{\Vect{Y}}{\Vect{X}_2} \right|.
\end{equation}
Now, notice that 
$$
\IM{\Vect{Y}}{\Vect{X}_1} \geq \IM{\Vect{Y}}{\Vect{X}_2}+\IM{\Vect{X}_1}{\Vect{X}_2}-\ent{\Vect{X}_2},
$$
from which one can deduce
\begin{equation} \label{eq-diffIM}
\left| \IM{\Vect{Y}}{\Vect{X}_1}-\IM{\Vect{Y}}{\Vect{X}_2} \right|  \leq \max\left( \ent{\Vect{X}_1},\ent{\Vect{X}_2}\right) -\IM{\Vect{X}_1}{\Vect{X}_2} = 
\max\left( \ent{\Given{\Vect{X}_1}{\Vect{X}_2}}, \ent{\Given{\Vect{X}_2}{\Vect{X}_1}}\right) = \DIB[I]{\Vect{X}_1}{\Vect{X}_2}.
\end{equation}
The result is then obtained by combining (\ref{hyp1Pen}), (\ref{eq1-proofPropDiff}) and~(\ref{eq-diffIM}).

$(ii)$ We can obtain the following result
\begin{eqnarray*}
\left| \NIB{\Vect{Y}}{\Vect{X}_1}- \NIB{\Vect{Y}}{\Vect{X}_2} \right| &\leq&  
\frac{\min(\pen{\Vect{Y}}{\Vect{X}_1}, \pen{\Vect{Y}}{\Vect{X}_2})\left(\left| \IM{\Vect{Y}}{\Vect{X}_1}- \IM{\Vect{Y}}{\Vect{X}_2} \right|+   \left| \pen{\Vect{Y}}{\Vect{X}_1}- \pen{\Vect{Y}}{\Vect{X}_2} \right|\right)}{\pen{\Vect{Y}}{\Vect{X}_1} \pen{\Vect{Y}}{\Vect{X}_2}} \\
&\leq & \frac{\left| \IM{\Vect{Y}}{\Vect{X}_1}- \IM{\Vect{Y}}{\Vect{X}_2} \right|+ \left| \pen{\Vect{Y}}{\Vect{X}_1}- \pen{\Vect{Y}}{\Vect{X}_2} \right|}{\max \left( \pen{\Vect{Y}}{\Vect{X}_1}, \pen{\Vect{Y}}{\Vect{X}_2} \right) }.
\end{eqnarray*}
The result then comes from (\ref{hyp1Pen}), (\ref{hyp2Pen}) and~(\ref{eq-diffIM}).
\end{proof}

Let us apply the previous result to our different examples:

\begin{corollary}

Let ${\Vect{X}_1},{\Vect{X}_2} \in \Gamma_{\Theta}$ with $\Theta=[c_1,c_2]$ and define $\gamma_i$ ($i=1,2$) such that $c_i=\gamma_i \ent{\Vect{Y}}$, then 
\begin{equation} \label{eq-resGenRedundant}
\left| \IB[\bullet]{\Vect{Y}}{\Vect{X}_1}- \IB[\bullet]{\Vect{Y}}{\Vect{X}_2} \right| \leq \left(1+\kappa_{1,\Theta}^\bullet\right)\; \DIB[I]{\Vect{X}_1}{\Vect{X}_2} 
\quad \mbox{ and } 
\left| \NIB[\bullet]{\Vect{Y}}{\Vect{X}_1}- \NIB[\bullet]{\Vect{Y}}{\Vect{X}_2} \right| \leq \frac{1+\kappa_{1,\Theta}^\bullet}{\kappa_{2,\Theta}^\bullet}\; \dIB[I]{\Vect{X}_1}{\Vect{X}_2} 
\end{equation}
where $\bullet$ stands for $S,R,P$ and $D$, and where
the different constants are expressed by 

\begin{center}
\begin{tabular}{|c|c|c|}
\hline
$\bullet$  & $\kappa_{1,\Theta}^\bullet$ & $\kappa_{2,\Theta}^\bullet$ \\
\hline
$S$  &   ${\alpha^{\vee}}$ & $(1-\alpha)+\alpha \gamma_{1,2}$ \\
\hline
$R$  &  ${\alpha^{\vee}}^2+ \frac{\alpha(1-\alpha)}{\sqrt{\gamma_1}}$ & $\left( (1-\alpha)+\alpha \sqrt{\gamma_{1,2}} \right)^2$ \\
\hline
$P$  & $\max\left( \frac{1-\alpha}{\gamma_1^{\alpha}},\frac{\alpha}{\gamma_1^{1-\alpha}},\mathbf{1}_{]0,1]}(\gamma_1)\right)$ & $\gamma_{1,2}^\alpha$\\
\hline
$D$  & $\frac{{\alpha^{\vee}}}{({\alpha^\wedge})^2} \; \frac1{(1+\gamma_{1,2})^2}$ & $\left(  \frac\alpha{\gamma_{1,2}} + (1-\alpha) \right)^{-1}$ \\
\hline
\end{tabular}
\end{center}
with $\gamma_{1,2}= \min \left(\gamma_1,\frac1{\gamma_2} \right)$.

\end{corollary}

\begin{proof}
For the sake of simplicity, let us denote by $m_i= \min\left(\ent{\Vect{Y}},\ent{\Vect{X}_i}\right)$ (resp. $m=\min\left(\ent{\Vect{Y}},\ent{\Vect{X}}\right) $) and by $M_i= \max\left(\ent{\Vect{Y}},\ent{\Vect{X}_i}\right)$ for $i=1,2$ (resp. $M=\max\left(\ent{\Vect{Y}},\ent{\Vect{X}}\right)$). Let us notice on the one hand that $|M_1+m_1-(M_2+m_2)|=|\ent{\Vect{X}_1}-\ent{\Vect{X}_2}|$ and on the other hand that 
$$
m \geq \left\{
\begin{array}{c}
\min\left( 1, \gamma_1\right)  \\
\min\left( 1, \frac1{\gamma_2}\right)
\end{array}
\right\} \geq \min \left(1,\gamma_1,\frac1{\gamma_2} \right) M = \gamma_{1,2} M
$$

\begin{list}{$\bullet$}{}
\item \Compl \penempty[S]: we have
\begin{eqnarray*}
| \pen[S,\alpha]{\Vect{Y}}{\Vect{X}_1}-\pen[S,\alpha]{\Vect{Y}}{\Vect{X}_2} |&=& |\alpha m_1+(1-\alpha)M_1-\alpha m_2 -(1-\alpha)M_2|\\
&=& |\alpha(m_1-m_2)+ (1-\alpha)(M_1-M_2)| \\
&\leq & {\alpha^{\vee}} | \ent{\Vect{X}_1} - \ent{\Vect{X}_2}|
\end{eqnarray*}
Moreover,
$$
\pen[S,\alpha]{\Vect{Y}}{\Vect{X}}= \alpha m + (1-\alpha)M \geq \left((1-\alpha)+\alpha \gamma_{1,2}\right) M
$$

\item \Compl \penempty[R]: we have
$$
| \pen[R]{\Vect{Y}}{\Vect{X}_1}-\pen[R]{\Vect{Y}}{\Vect{X}_2} | = \left| \alpha^2 (m_1-m_2) + (1-\alpha)^2 (M_1-M_2) + 2\alpha(1-\alpha) \sqrt{\ent{\Vect{Y}}} \left(\sqrt{\ent{\Vect{X}_1}}- \sqrt{\ent{\Vect{X}_2}}\right)
\right| 
$$
Furthermore, we may obtain
$$\left| \alpha^2 (m_1-m_2) + (1-\alpha)^2 (M_1-M_2)\right|  \leq {\alpha^{\vee}}^2|\ent{\Vect{X}_1}-\ent{\Vect{X}_2}|
$$
and
$$
\left| \sqrt{\ent{\Vect{Y}}} (\sqrt{\ent{\Vect{X}_1}}- \sqrt{\ent{\Vect{X}_2}}) \right| = \frac{\sqrt{\ent{\Vect{Y}}}}{2 \sqrt{\min(\ent{\Vect{X}_1},\ent{\Vect{X}_2})}}\times |\ent{\Vect{X}_1}-\ent{\Vect{X}_2}| \leq \frac1{2\sqrt{\gamma_1}} |\ent{\Vect{X}_1}-\ent{\Vect{X}_2}|.
$$
Hence,
$$
| \pen[R]{\Vect{Y}}{\Vect{X}_1}-\pen[R]{\Vect{Y}}{\Vect{X}_2} |\leq \left( {\alpha^{\vee}}^2  + \frac{\alpha(1-\alpha)}{\sqrt{\gamma_1}} \right) |\ent{\Vect{X}_1}-\ent{\Vect{X}_2}|.
$$

Moreover, one can prove 
$$
\pen[R,\alpha]{\Vect{Y}}{\Vect{X}}= (\alpha \sqrt{m} + (1-\alpha)\sqrt{M})^2 \geq \left((1-\alpha)+\alpha \sqrt{\gamma_{1,2}} \right)^2 M
$$
\item \Compl \penempty[P]: we have (by assuming $\ent{\Vect{X}_2}>\ent{\Vect{X}_1}$)
\begin{eqnarray*}
| \pen[P,\alpha]{\Vect{Y}}{\Vect{X}_1}-\pen[P,\alpha]{\Vect{Y}}{\Vect{X}_2} | &=& \left| m_1^\alpha M_1^{1-\alpha} - m_2^\alpha M_2^{1-\alpha} \right| \\
&=& \left\{ 
\begin{array}{ll}
\ent{\Vect{Y}}^\alpha \left( \ent{\Vect{X}_2}^{1-\alpha}-\ent{\Vect{X}_1}^{1-\alpha} \right) & \mbox{ if } \ent{\Vect{Y}}\leq \min\left(\ent{\Vect{X}_1},\ent{\Vect{X}_2}\right) \\
\ent{\Vect{Y}}^{1-\alpha} \left( \ent{\Vect{X}_2}^{\alpha}-\ent{\Vect{X}_1}^{\alpha} \right) & \mbox{ if } \ent{\Vect{Y}}\geq \max\left(\ent{\Vect{X}_1},\ent{\Vect{X}_2}\right) \\
\ent{\Vect{Y}}^\alpha  \ent{\Vect{X}_2}^{1-\alpha}-\ent{\Vect{X}_1}^{\alpha}\ent{\Vect{Y}}^{1-\alpha} & \mbox{otherwise.}
\end{array}
\right. 
\end{eqnarray*}
Note that the third case cannot occur if $\gamma_1\geq 1$.
\begin{eqnarray*}
| \pen[P,\alpha]{\Vect{Y}}{\Vect{X}_1}-\pen[P,\alpha]{\Vect{Y}}{\Vect{X}_2} |
&\leq& \left\{
\begin{array}{ll}
\frac{1-\alpha}{\gamma_1^{\alpha}} (\ent{\Vect{X}_2}-\ent{\Vect{X}_1})& \mbox{ if } \ent{\Vect{Y}}\leq \min\left(\ent{\Vect{X}_1},\ent{\Vect{X}_2}\right) \\
\frac{\alpha}{\gamma_1^{1-\alpha}}(\ent{\Vect{X}_2}-\ent{\Vect{X}_1}) & \mbox{ if } \ent{\Vect{Y}}\geq \max\left(\ent{\Vect{X}_1},\ent{\Vect{X}_2}\right) \\
\ent{\Vect{X}_2}-\ent{\Vect{X}_1} & \mbox{otherwise}
\end{array} \right. \\
&\leq& \max\left( \frac{1-\alpha}{\gamma_1^{\alpha}},\frac{\alpha}{\gamma_1^{1-\alpha}},\mathbf{1}_{]0,1]}(\gamma_1)\right) |\ent{\Vect{X}_2}-\ent{\Vect{X}_1}|. 
\end{eqnarray*}
Moreover, we may obtain 
$$
\pen[P,\alpha]{\Vect{Y}}{\Vect{X}} =m^\alpha M^{1-\alpha} \geq \gamma_{1,2}^\alpha M
$$

\item \Compl \penempty[D]: we have
\begin{eqnarray}
| \pen[D,\alpha]{\Vect{Y}}{\Vect{X}_1}-\pen[D,\alpha]{\Vect{Y}}{\Vect{X}_2} | & =& \frac{\alpha M_1 M_2 (m_1-m_2) + (1-\alpha) m_1 m_2 (M_1-M_2)}{(\alpha M_1 + (1-\alpha)m_1)(\alpha M_2 + (1-\alpha)m_2)} \\
&\leq & \frac{{\alpha^{\vee}}}{({\alpha^\wedge})^2} \frac{M_1 M_2}{(m_1+M_1)(m_2+M_2)} \left| \ent{\Vect{X}_2}-\ent{\Vect{X}_1}\right| \\
&\leq& \frac{{\alpha^{\vee}}}{({\alpha^\wedge})^2} \; \frac1{(1+\gamma_{1,2})^2} \left| \ent{\Vect{X}_2}-\ent{\Vect{X}_1}\right|
\end{eqnarray}
Finally, we also have 
$$
\pen[D,\alpha]{\Vect{Y}}{\Vect{X}} = \left(  \frac\alpha{m} + \frac{1-\alpha}{M} \right)^{-1}  \geq
 \left(  \frac\alpha{\gamma_{1,2}} + (1-\alpha) \right)^{-1} M. 
$$
\end{list}
\end{proof}

\begin{remark}
Note that when $\alpha\leq \frac12$, the measure $\IBempty[S,\alpha]$ is a metric and so we derive~(\ref{eq-resGenRedundant}) directly from \textbf{[P3bis]}.
\end{remark}


\bibliographystyle{plainnat.bst}
\bibliography{crit}

\newpage
\noindent \textbf{Authors}: Jean-Fran\c{c}ois Coeurjolly, Rémy Drouilhet and Jean-François Robineau. \\
\noindent \textbf{Address}: LABSAD, BSHM, 1251 avenue centrale BP 47 - 38040 GRENOBLE Cedex 09. \\
\noindent \textbf{E-mail addresses}:
\begin{itemize}
\item[] \texttt{Jean-Francois.Coeurjolly@upmf-grenoble.fr}
\item[] \texttt{Remy.Drouilhet@upmf-grenoble.fr}
\end{itemize}
\noindent \textbf{Corresponding author}: Jean-François Coeujolly.

\end{document}